\documentstyle[12pt]{amsart}

\def\commdiag#1{\null\,
  \vcenter{\commdiagbaselines
  \m@th\ialign{\hfil$##$\hfil&&\hfil$\mkern4mu ##$\hfil\crcr
      \mathstrut\crcr\noalign{\kern-\baselineskip}
      #1\crcr\mathstrut\crcr\noalign{\kern-\baselineskip}}}\,}%

\def\commdiagbaselines{\baselineskip15pt \lineskip3pt \lineskiplimit3pt }%
\def\gridcommdiag#1{\null\,
  \vcenter{\offinterlineskip
  \m@th\ialign{&\vbox to\vgrid{\vss
    \hbox to\hgrid{\hss\smash@@{$##$}\hss}}\crcr
      \mathstrut\crcr\noalign{\kern-\vgrid}
      #1\crcr\mathstrut\crcr\noalign{\kern-.5\vgrid}}}\,}%

\newdimen\harrowlength \harrowlength=60pt
\newdimen\varrowlength \varrowlength=.618\harrowlength
\newdimen\sarrowlength \sarrowlength=\harrowlength

\newdimen\hgrid \hgrid=15pt
\newdimen\vgrid \vgrid=15pt

\newdimen\hchannel  \hchannel=0pt
\newdimen\vchannel  \vchannel=0pt
\newdimen\channelwidth \channelwidth=3pt

\input epsf
\theoremstyle{plain}
\newtheorem{thm}[subsection]{Theorem}
\newtheorem{lem}[subsection]{Lemma}
\newtheorem{prop}[subsection]{Proposition}
\newtheorem{cor}[subsection]{Corollary}

\theoremstyle{definition}
\newtheorem{rem}[subsection]{Remark}
\newtheorem{defn}[subsection]{Definition}

\def\v{\vee}
\def\w{\wedge}

\def\L{\cal{L}}

\def\e{\varepsilon}
\def\a{\alpha}
\def\t{\tau}

\def\f{\phi}
\def\l{\lambda}

\def\h{\chi}
\def\m{\mu}

\def\v{\vee}
\def\w{\wedge}

\def\L{\cal{L}}

\def\e{\varepsilon}
\def\a{\alpha}
\def\t{\tau}

\def\f{\phi}
\def\l{\lambda}

\def\h{\chi}
\def\m{\mu}

\def\vs{\vskip}
\def\ni{\noindent}

\begin{document}
\title[Singular loci of certain toric varieties]{Singular loci of Hibi toric varieties}
\address{Department of Mathematics\\ Northeastern University\\ Boston, MA 02115}

\author[V. Lakshmibai]{V. Lakshmibai${}^{\dag}$}
\author{H. Mukherjee}
\thanks{${}^{\dag}$ Partially suported
by NSF grant DMS-0652386 and Northeastern University RSDF 07-08.}
\maketitle
\begin{abstract}
We first construct explicit bases for the cotangent spaces at
singular points on Hibi toric varieties, i.e., toric varieties
associated to finite distributive lattices. We then determine the
singular loci of  these toric varieties.
\end{abstract}

\section*{Introduction}
Let $K$ denote the base field which we assume to be algebraically
closed of arbitrary characteristic. Given a distributive lattice
$\mathcal{L}$, let $X_{\mathcal{L}}$ denote the \emph{Hibi
variety} - the affine variety in $\mathbb{A}^{\#\mathcal{L}}$
whose vanishing ideal is generated by the binomials $X_\tau
X_\varphi -X_{\tau\vee\varphi}X_{\tau\w\varphi}$ in the polynomial
algebra $K[X_\alpha,\alpha\in\mathcal{L}]$; here,
${\tau\vee\varphi}$ (resp. ${\tau\w\varphi}$) denotes the
\emph{join} - the smallest element of $\mathcal{L}$ greater than
both $\tau,\varphi$ (resp. the \emph{meet} - the largest element
of $\mathcal{L}$ smaller than both $\tau,\varphi$). In the sequel,
we shall refer to the quadruple $(\t,\varphi,{\tau\vee\varphi},
{\tau\w\varphi})$ ($(\t,\varphi)$ being a skew pair) as a
\emph{diamond}. These varieties were extensively studied by Hibi
in \cite{Hi} where Hibi proves that $X_{\mathcal{L}}$ is a normal
variety. On the other hand, Eisenbud-Sturmfels show in \cite{ES}
that a binomial prime ideal is toric (here, ``toric ideal" is in
the sense of \cite{St}). Thus one obtains that $X_{\mathcal{L}}$
is a normal toric variety. We refer to such an $X_{\mathcal{L}}$
as a Hibi toric variety.

For $\mathcal{L}$ being the Bruhat poset of Schubert varieties in
a minuscule $G/P$, it is shown in \cite{GLdef} that
$X_{\mathcal{L}}$ flatly deforms to $G/P$, i.e., there exists a
flat family over ${\mathbb{A}}^{1}$ with $G/P$ as the generic
fiber and $X_{\mathcal{L}}$ as the special fiber. More generally
for a Schubert variety $X(w)$ in a minuscule $G/P$, it is shown in
\cite{GLdef} that $X_{{\mathcal{L}}_w}$ flatly deforms to $X(w)$
(here, ${\mathcal{L}}_w$ is the Bruhat poset of Schubert
subvarieties of $X(w)$). In a subsequent paper (cf. \cite{g-l}),
the authors of loc.cit., studied the singularities of
$X_{\mathcal{L}}, \mathcal{L}$ being the Bruhat poset of Schubert
varieties in the Grassmannian; further, in loc.cit., the authors
gave a conjecture giving a necessary and sufficient condition for
a point on $X_{\mathcal{L}}$ to be smooth, and proved in loc.cit.,
the sufficiency part of the conjecture. Subsequently, the
necessary part of the conjecture was proved in \cite{Font1} by
Batyrev et al. The toric varieties $X_{\mathcal{L}}, \mathcal{L}$
being the Bruhat poset of Schubert varieties in the Grassmannian
play an important role in the area of mirror symmetry; for more
details, see \cite{Font2,Font1}. We refer to such an
$X_{\mathcal{L}}$ as a \emph{Grassmann-Hibi}  toric variety.

In this paper, we first determine an explicit description (cf.
\S\ref{cone}) of the cone $\sigma$ associated to a Hibi toric
variety and the dual cone $\sigma^{\vee}$ (note that if
$S_{\sigma}$ denotes the semigroup of integral points in
$\sigma^{\vee}$, then $K[X_{\mathcal{L}}]$, the $K$-algebra of
regular functions on $X_{\mathcal{L}}$ can be identified with the
semigroup algebra $K[S_{\sigma}]$). Using this we determine the
cotangent space at any point on $X_{\mathcal{L}}$ very explicitly
as described below.

For each face $\tau$ of $\sigma$, there exists a distinguished
point $P_\tau $ in the torus orbit ${\mathcal{O}}_\tau$
corresponding to $\tau$; namely, identifying a (closed) point in
$X_{\mathcal{L}}$ with a semigroup map $S_{\sigma}\rightarrow
K^*\cup\{0\}$,  $P_\tau $ corresponds to the semigroup map
$f_\tau$ which sends $u\in S_{\sigma}$ to $1$ or $0$ according as
$u$ is in $\tau^{\perp}$ or not. We have that $f_\tau$ induces an
algebra map $\psi_\tau:K[X_{{\mathcal{L}}}]\rightarrow K$ whose
kernel is precisely $M_{P_{\tau}}$, the maximal ideal
corresponding to $P_{\tau}$. For $\alpha\in{\mathcal{L}}$, we
shall denote by $P_{\tau}(\alpha)$, the co-ordinate of $P_{\tau}$
(considered as a point of $\mathbb{A}^{\#\mathcal{L}}$)
corresponding to $\alpha$. Set
$$D_\tau=\{\alpha\in{\mathcal{L}}\,|\,P_{\tau}(\alpha)\not= 0\}$$
It turns out that $D_\tau$ is an embedded sublattice of
${\mathcal{L}}$; further, $D_{\tau}$ determines the local behavior
at $P_{\tau}$. To make this more precise, identifying
$K[X_{\mathcal{L}}]$ as the quotient of the polynomial algebra
$K[X_\alpha,\alpha\in\mathcal{L}]$ by the ideal generated by the
binomials $X_\tau X_\varphi
-X_{\tau\vee\varphi}X_{\tau\w\varphi}$, let $x_\alpha$ denote the
image of $X_\alpha$ in $K[X_{\mathcal{L}}]$. For $\alpha\in
\mathcal{L}$, set
$$F_\alpha=\begin{cases}x_\alpha,&{\mathrm{\ if\ }}
{\alpha}\not\in D_{\tau}\\
1-x_\alpha,&{\mathrm{\ if\ }}\alpha\in D_{\tau}
\end{cases}$$ Then we have that $M_{P_{\tau}}$ is generated by
$\{F_\alpha,\alpha\in{\mathcal{L}}\}$. Fix a maximal chain
$\Gamma$ in $D_\tau$. Let $\Lambda_\tau(\Gamma)$ denote the
sublattice of ${\mathcal{L}}$ consisting of all maximal chains of
${\mathcal{L}}$ containing $\Gamma$. Let $E_\tau$ denote the set
of all $\alpha\in{\mathcal{L}}$ such that there exists a $\beta\in
D_\tau$ such that $(\alpha,\beta)$ is a diagonal of a diamond
whose other diagonal is contained in $
{\mathcal{L}}\,\setminus\,D_\tau$. Let
$Y_\tau(\Gamma)=\Lambda_\tau(\Gamma)\,\setminus\,E_\tau$. For
$F\in M_{P_{\tau}}$, let ${\overline{F}}$ denote the class of $F$
in $M_{P_{\tau}}\,/\,M_{P_{\tau}}^2$. We define an equivalence
relation (cf. \S \ref{equ}) on ${\L}\,\setminus\,D_\tau$ in such a
way that for two elements $\theta,\theta'$ in the same equivalence
class, we have ${\overline{F}}_{\theta}=
{\overline{F}}_{\theta'}$. Denoting by $G_\tau(\Gamma)$ the set of
equivalence classes $[\theta]$, where $\theta\in
Y_\tau(\Gamma)\,\setminus\,\{\Gamma\cup E_\tau\}$ and
${\overline{F}}_{[\theta]}$ is non-zero in
$M_{P_{\tau}}\,/\,M_{P_{\tau}}^2$, we have

\vs.1cm\ni{\textbf{Theorem 1:}} (cf. Theorem \ref{main2})
$\{{\overline{F}}_{[\theta]},[\theta]\in
G_{\tau}(\Gamma)\}\dot\cup\{{\overline{F}}_\gamma,\gamma\in\Gamma\}$
is a basis for the cotangent space
$M_{P_{\tau}}\,/\,M_{P_{\tau}}^2$.

Using the above result, we determine Sing$\,X_{\L}$ (cf. Theorem
\ref{main3}):

\ni Let
${\mathcal{S_{\L}}}=\{\tau<\sigma\,|\,\#G_{\tau}(\Gamma)+\#\Gamma>\#{\L}\}$.

\vs.1cm\ni{\textbf{Theorem 2:}} (cf. Theorem \ref{main3})
Sing$\,X_{\L}={\underset{\tau\in\mathcal{S_{\L}}}{\cup}}\,X(D_\tau)$.

\ni (Here, $X(D_\tau)$ denotes the Hibi variety associated to the
distributive lattice $D_\tau$.)

\vs.2cm\ni{\textbf{Sketch of the proof of the main theorem:}}
Using our explicit description of the generators for
$\sigma,\sigma^{\vee}$ (cf. \S \ref{cone}), we first determine
explicitly the embedded sublattice associated to a face $\tau$ of
$\sigma$. We then analyze the local expression around $P_{\tau}$
for any $f\in I(X_{{\mathcal{L}}})$, the vanishing ideal of
$X_{\mathcal{L}}$, and show the generation of the degree one part
of gr$\,(R_{{\mathcal{L}}}, M_{P_{\tau}})$ by
$\{{\overline{F}}_\theta,\theta\in Y_\tau(\Gamma)\}$ (here,
$R_{{\mathcal{L}}}=K[X_{{\mathcal{L}}}]$, and $M_{P_{\tau}}$ is
the maximal ideal in $R_{{\mathcal{L}}}$ corresponding to
$P_{\tau}$). We then define the equivalence relation on
${\L}\,\setminus\,D_\tau$. The linear independence of
$\{{\overline{F}}_{[\theta]}\,[\theta]\in G_\tau(\Gamma)\}\cup
\{{\overline{F}}_\gamma,\gamma\in\Gamma \}$ in
$M_{P_{\tau}}\,/\,M_{P_{\tau}}^2$ is proved using the defining
equations of $X_{\mathcal{L}}$ (as a closed subvariety of
${\mathbb{A}}^{\# {\mathcal{L}}}$), thus proving Theorem 1.
Theorem 2 is then deduced from Theorem 1.

The sections are organized as follows: In \S \ref{sec1}, we recall
generalities on toric varieties. In \S \ref{varlat}, we recall some
basic results on distributive lattices. In \S \ref{5.13}, we
introduce the Hibi toric variety $X_{{\mathcal{L}}}$, and recall
some basic results on $X_{{\mathcal{L}}}$. In \S \ref{cone}, we
determine generators for the cone $\sigma$ (and the dual cone
$\sigma^{\vee}$); further, for a face $\tau$ of $\sigma$, we
introduce $D_\tau$ and derive some properties of $D_\tau$. In \S
\ref{tan-cot}, we prove the generation of
$M_{P_{\tau}}\,/\,M_{P_{\tau}}^2$ by
$\{{\overline{F}}_{[\theta]},[\theta]\in
G_{\tau}(\Gamma)\}\dot\cup\{{\overline{F}}_\gamma,\gamma\in\Gamma\}$,
while the linear independence is proved in \S \ref{cotgt}.

\section{Generalities on toric varieties}\label{sec1} Since our main object of
study is a certain affine toric variety, we recall in this section
some basic definitions on affine toric varieties. Let $T=(K^*)^m$
be an $m$-dimensional torus.
\begin{defn}\label{equi} (cf. \cite{F}, \cite{toroidal})
An {\em equivariant affine embedding}  of a torus $T$ (or also
\emph{equivariant affine toroidal embedding}) is an affine variety
$X\subseteq{\mathbb{A}}^l$ containing $T$ as a dense open subset
and equipped with a $T$-action $T\times X\to X$ extending the
action $T\times T\to T$ given by multiplication. If in addition
$X$ is normal, then $X$ is called an {\em affine toric variety}.
\end{defn}

\subsection{Resum\'e of combinatorics of affine toric
varieties}\label{comb} Let $M$ be the character group of $T$, and
$N$ the $\mathbb{Z}$-dual of $M$. Let $X$ be an equivariant affine
toroidal embedding of $T$; let $R=K[X]$ (the ring of regular
functions on $X$). We note the following:

$\bullet$ \textbf{$M$-gradation:} We have an action of $T$ on $R$,
and hence writing $R$ as a sum of $T$ weight spaces, we obtain a
$M$-grading for $R$: $R={\underset{\{\chi\in
M\}}{\oplus}}R_{\chi}$, where

\ni $$R_{\chi}=\{f\in R|tf=\chi(t)f, \forall t\in T\}$$

$\bullet$ \textbf{The semi group $S$:} Let $S=\{\chi\in
M|R_{\chi}\neq 0\}$. Then via the multiplication in $R$, $S$
acquires a semi group structure. Thus $S$ is a sub semigroup of
$M$, and $R$ gets identified with the semi group algebra $K[S]$.

$\bullet$ \textbf{Finite generation of $S$:} In view of the fact
that $R$ is a finitely generated $K$-algebra, we obtain that $S$
is a finitely generated sub semigroup of M.

$\bullet$ \textbf{Generation of $M$ by $S$:} The fact that $T$ and
$X$ have the same function field (since $T$ is a dense open subset
of $X$) implies that $M$ is generated by $S$ (note that $K[T]$ is
the group algebra $K[M]$).

Thus given  an equivariant affine toroidal embedding  $X$ (of the
torus $T$), $X$ determines a finitely generated,
 sub semigroup $S$ of $M$ which generates $M$ in such a
way that the semi group algebra $K[S]$ is the ring of regular
functions on $X$. Conversely, given a finitely generated, sub
semigroup $S$ of $M$ which generates $M$, we have clearly that
$X:=Spec\,K[S]$ is an equivariant affine toroidal embedding of
$T$.

Thus we obtain a bijection between \{equivariant affine toroidal
embeddings of $T$\} and \{finitely generated,
 sub semigroups of $M$ which generate $M$\}.

 Let notation be as above.

$\bullet$ \textbf{Saturation of $S$:} We have (see \cite{toroidal}
for example) that the ring $R$ (being identified as the semi group
algebra $K[S]$) is normal if and only if $S$ is saturated (recall
that a sub semigroup $A$ in $M$ is \emph{saturated},  if for $\chi
\in M $, $r\chi \in A\ \Longrightarrow \ \chi \in A$ (where $r$ is
an integer $>1$; equivalently, $A$ is precisely the lattice points
in the cone generated by $A$, i.e., $A=\theta\cap A$ where
$\theta=\sum a_ix_i, a_i \in{\mathbb{R}}_+,x_i\in A$)).

Thus given an affine toric variety $X$ (with torus action by $T$),
$X$ determines a finitely generated, saturated sub semigroup $S$
of $M$ which generates $M$ in such a way that the semi group
algebra $K[S]$ is the ring of regular functions on $X$.

Let $S$ be as above, namely, a finitely generated, saturated, sub
semigroup of $M$ which generates $M$. Let $\theta:=\sum a_ix_i,
a_i \in{\mathbb{R}}_+,x_i\in S$ be the cone generated by $S$. Then
$\theta$ is a \emph{rational convex polyhedral cone} (i.e., a cone
with a (finite) set of lattice points as generators); further,
$\theta$ is not contained in any hyperplane in
$M_{\mathbb{R}}:=M\otimes\mathbb{R}$ (since $S$ generates $M$, we
have that $\theta$ generates $M\otimes{\mathbb{R}}$). Let
$$\sigma:=\theta^{\v}=\{v\in N_{\mathbb{R}}|v(f)\ge 0,\forall f\in \theta\}$$
(note that $N_{\mathbb{R}} (=N\otimes{\mathbb{R}})$ is the linear
dual of $M_{\mathbb{R}}$). Then $\sigma$ is a \emph{strongly
convex rational polyhedral cone} (one may take this to be the
definition of a strongly convex rational convex polyhedral cone,
namely, a rational convex polyhedral cone $\tau$ is \emph{strongly
convex} if $\tau^{\v}$ is not contained in any hyperplane;
equivalently, $\tau$ does not contain any non-zero linear
subspace.)

Thus an affine toric variety $X$ (with a torus action by $T$)
determines a strongly convex, rational polyhedral cone $\sigma$
(inside $N_{\mathbb{R}}(=M_{\mathbb{R}}^{*})$) in such a way that
$K[X]$ is the semigroup algebra $K[S_\sigma]$, where $S_\sigma$ is
the sub semigroup in $M$ consisting of the set of lattice points
in $\sigma^{\v}$.

Conversely, starting with a strongly convex rational polyhedral
cone $\sigma$,

\ni $S_\sigma:=\sigma^{\v}\cap M$ is a finitely generated,
saturated, sub semigroup of $M$ which generates $M$; and hence
determines an affine toric variety $X$ (with a torus action by $T$),
namely, $X=Spec\,K[S_\sigma]$.

This sets up a bijection between \{affine toric varieties with torus
action by $T$\}, \{finitely generated, saturated, sub semigroups of
$M$ which generate $M$\}, and \{strongly convex rational polyhedral
cones inside $N_{\mathbb{R}}$\}.

In the sequel, we shall denote $Spec\,K[S_\sigma]$ by $X_\sigma$.


\subsection{Toric ideals, binomial prime ideals and equivariant
affine toroidal embeddings:}\label{toric}
 In this subsection we shall establish bijections
between \{equivariant affine toroidal embeddings\}, \{toric
ideals\} and \{binomial prime ideals\} (here, by a \emph{a
binomial ideal} we mean an ideal generated by binomials, i.e.,
polynomials with at most two terms).

Let ${\mathcal{A}}=\{\chi_1,\dots,\chi_l\}$ be a subset of
${\mathbb{Z}}^d$. Consider the  map
\begin{equation*}
\pi_{\mathcal{A}}:{\mathbb{Z}}_+^l\to {\mathbb{Z}}^d,\quad
{\mathbf{u}}=(u_1,\dots,u_l)\mapsto u_1\chi_1+\dots+u_l\chi_l.
\end{equation*}
Let $K[{\mathbf{x}}]:=K[x_1,\dots,x_l]$, $K[{\mathbf{t}}^{\pm
1}]:=K[t_1,\dots,t_d, t_1^{-1},\dots,t_d^{-1}]$.

 The map $\pi_{{\mathcal{A}}}$ induces a homomorphism of
semigroup algebras
\begin{equation*}
\widehat{\pi_{{\mathcal{A}}}}:K[{\mathbf{x}}]\to
K[{\mathbf{t}}^{\pm 1}],\qquad x_i\mapsto{\mathbf{t}}^ {\chi_i}.
\end{equation*}

\ni (If $\chi_i=(a_{i1},\cdots,a_{id})$, then ${\mathbf{t}}^
{\chi_i}=t_1^{a_{i1}}\cdots t_d^{a_{id}}$.)

\begin{defn}\label{tor.ideal}(cf. \cite{St}) The kernel of $\widehat{\pi}_{{\mathcal{A}}}$ is denoted by
$I_{{\mathcal{A}}}$  and called the\/ {\em toric ideal}
 associated to $\mathcal{A}$.
\end{defn}

Note that a toric ideal is prime. We shall now show that we have a
bijection between \{toric ideals\} and \{equivariant affine
toroidal embeddings\}

 Consider the action of $T(=(K^*)^d)$ on ${{\mathbb{A}}}^l$ given by
$${\mathbf{t}}(a_1,\cdots,\a_l)=({\mathbf{t}}^{\chi_1}a_1,\cdots,{\mathbf{t}}^{\chi_l}a_l)$$
 Then
${\mathcal{V}}(I_{{\mathcal{A}}})$, the affine variety of the zeroes
in $K^l$ of $I_{{\mathcal{A}}}$, is simply the Zariski closure of
the $T$-orbit through $(1,1,\dots ,1)$. Let $M_{{\mathcal{A}}}$ be
the ${\mathbb{Z}}$-span of ${\mathcal{A}}$ (inside
$X(T)={{\mathbb{Z}}}^d$), and $d_{{\mathcal{A}}}$ the rank of
$M_{{\mathcal{A}}}$. Let $S_{{\mathcal{A}}}$ be the sub semigroup of
${{\mathbb{Z}}}^d$ generated by ${\mathcal{A}}$. Let
$T_{{\mathcal{A}}}$ be the $d_{{\mathcal{A}}}$-dimensional torus
with $M_{{\mathcal{A}}}$ as the character group. Then we have

\begin{prop}\label{zar} ${\mathcal{V}}(I_{{\mathcal{A}}})$ is an
equivariant affine embedding of $T_{{\mathcal{A}}}$ (of dimension
$d_{{\mathcal{A}}}$); further, $K[{\mathcal{V}}(I_{{\mathcal{A}}})]$
is the semigroup algebra $K[S_{{\mathcal{A}}}]$.
\end{prop}

\begin{rem} In the above definition, we do not require ${\mathcal{V}}(I_{{\mathcal{A}}})$
to be normal. Note that ${\mathcal{V}}(I_{{\mathcal{A}}})$ is normal
if and only if $S_{{\mathcal{A}}}$ is saturated.  A variety of the
form ${\mathcal{V}}(I_{{\mathcal{A}}})$, is also sometimes called
an\/ {\em affine toric variety}. But in this paper, by an affine
variety, we shall mean a normal toric variety as defined in
Definition \ref{equi}.
\end{rem}

\subsubsection{The vanishing ideal of an equivariant affine toroidal embedding:}\label{vanishing}
Conversely, let $X$ be an equivariant affine toroidal embedding
(with a torus action  by $T(=(K)^*)^d$ as in Definition \ref{equi}).
Let $K[X]=K[S]$, for some suitable finitely generated sub semigroup
of $X(T)(={\mathbb{Z}}^{d})$. Fix a set of generators
${\mathcal{A}}:=\{\chi_1,\cdots,\chi_l\}$ for the sub semigroup $S$
of ${\mathbb{Z}}^{d}$. We have a surjection
$$K[x_1,\cdots,x_l]\rightarrow K[X],\ x_i\mapsto{\mathbf{t}^{\chi_i}} $$
whose kernel equals $I_{{\mathcal{A}}}$. Thus we obtain
\begin{prop}\label{above'} With $X$ as above, we have that
$I(X)$, the vanishing ideal of $X$ (for the embedding
$X\hookrightarrow{\mathbb{A}}^l$) is a toric ideal.
\end{prop}
Thus we obtain a bijection between \{equivariant, affine toroidal
embeddings\} and \{toric ideals\}.
Next we want to show that we have a bijection between \{equivariant,
affine toroidal embeddings\} and \{binomial prime ideals\}.

First, we recall the following (see \cite{St}).
\begin{prop} \label{vegen}
The toric ideal $I_{\mathcal{A}}$ is spanned as a  $K$-vector
space by the set of binomials
\begin{equation*}
\{{\mathbf{x}}^{\mathbf{u}}-\mathbf{x}^{\mathbf{v}}\mid
{\mathbf{u}},{\mathbf{v}}\in {\mathbb{Z}}_+^l \text{ with }
\pi_{\mathcal{A}}
(\mathbf{u})=\pi_{\mathcal{A}}(\mathbf{v})\}.\tag{*}
\end{equation*}
\end{prop}
As an immediate consequence, we obtain
\begin{cor}\label{toric'}
A toric ideal is a binomial prime ideal.
\end{cor}
In view of the above bijection, we have the following:

\ni\textbf{Reformulation of Corollary} \ref{toric'}: The vanishing
ideal of an equivariant, affine toroidal embedding is binomial
prime.

Conversely, we shall now show that a binomial prime ideal defines an
equivariant affine toroidal embedding.

\subsection{Varieties defined by binomials}\label{varbino}

Let $l\ge 1$, and let $X$ be an affine variety in
${\mathbb{A}}^l$, not contained in any of the coordinate
hyperplanes $\{ x_i=0\}$. Further, let $X$ be irreducible, and let
its defining prime ideal $I(X)$ be generated by $m$ binomials
\begin{equation*}
x_1^{a_{i1}}\dots x_l^{a_{il}}-\l_i x_1^{b_{i1}}\dots
x_l^{b_{il}}\, ,\qquad 1\le i\le m\quad .\tag{$*$}
\end{equation*}
Consider the natural action of the torus $T_l=(K^*)^l$ on
${\mathbb{A}}^l$,
\begin{equation*}
(t_1,\dots,t_l)\cdot (a_1,\dots,a_l)=(t_1a_1,\dots,t_la_l).
\end{equation*}
Let $X(T_l)=\text{Hom} (T_l,{\mathbb{G}}_m)$ be the character
group of $T_l$, and let $\e_i\in X(T_l)$ be the character
\begin{equation*}
\e_i(t_1,\dots,t_l)=t_i\, ,\qquad 1\le i\le l.
\end{equation*}
For $1\le i\le m$, let
\begin{equation*}
\varphi_i=\sum_{r=1}^l(a_{ir}- b_{ir})\e_r\, .
\end{equation*}
Set $T=\cap _{i=1}^m\ker \varphi_i$, and
$X^{\circ}=\{(x_1,\dots,x_l)\in X \mid x_i\ne 0, \text{ for all }
1\le i\le l\}$.

\begin{prop} \label{binotoric}
Let notations be as above.

$(1)$ There is a canonical action of $T$ on $X$.

$(2)$  $X^{\circ}$ is $T$-stable. Further, the action of $T$ on
$X^{\circ}$ is simple and transitive.

$(3)$ $T$ is a subtorus of $T_l$, and $X$ is an equivariant affine
embedding of $T$.
\end{prop}
\begin{pf}
$(1)$ We consider  the  (obvious) action of $T$ on
${\mathbb{A}}^l$. Let $(x_1,\dots,x_l)\in X$,
${\bold{t}}=(t_1,\dots,t_l)\in T$, and $(y_1,\dots,y_l)=
{\bold{t}}\cdot (x_1,\dots,x_l)=(t_1x_1,\dots,t_lx_l)$. Using the
fact that $(x_1,\dots,x_l)$ satisfies $(*)$, we obtain
\begin{equation*}
y_1^{a_{i1}}\dots y_l^{a_{il}}=t_1^{a_{i1}}\dots
t_l^{a_{il}}x_1^{a_{i1}}\dots x_l^{a_{il}}= \l_i t_1^{b_{i1}}\dots
t_l^{b_{il}}x_1^{b_{i1}}\dots x_l^{b_{il}}=\l_i y_1^{b_{i1}}\dots
y_l^{b_{il}},
\end{equation*}
for all $1\le i\le m$, i.e., $(y_1,\dots,y_l)\in X$. Hence
${\bold{t}}\cdot (a_1,\dots,a_l)\in X$ for all ${\bold{t}}\in T$.

$(2)$ Let ${\bold{x}}=(x_1,\dots,x_l)\in X^{\circ}$, and
${\bold{t}}=(t_1,\dots,t_l)\in T$. Then, clearly ${\bold{t}}\cdot
{\bold{x}}\in X^{\circ}$. Considering ${\bold{x}}$ as a point in
${\mathbb{A}}^l$, the isotropy subgroup in $T_l$ at ${\bold{x}}$
is $\{\text{id}\}$. Hence the isotropy subgroup in $T$ at $x$ is
also $\{\text{id}\}$. Thus the action of $T$ on $X^{\circ}$ is
simple.

Let $(x_1,\dots,x_l)$, $(x'_1,\dots,x'_l)\in X^{\circ}$. Set
${\bold{t}}=(t_1,\dots,t_l)$, where $t_i=x_i/x'_i$. Then, clearly
${\bold{t}}\in T$. Thus $(x_1,\dots,x_l)={\bold{t}}\cdot
(x'_1,\dots,x'_l)$. Hence the action of $T$ on $X^{\circ}$ is
simple and transitive.

$(3)$ Now, fixing a point ${\bold{x}}\in X^{\circ}$, we obtain
from $(2)$ that the orbit map ${\bold{t}}\mapsto {\bold{t}}\cdot
{\bold{x}}$ is in fact an isomorphism of $T$ onto $X^{\circ}$.
Also, since $X$ is not contained in any of the coordinate
hyperplanes, the open set $X_i=\{(x_1,\dots,x_l)\in X\mid x_i\ne
0\}$ is nonempty for all $1\le i\le l$. The irreducibility of $X$
implies that the sets $X_i$, $1\le i\le l$, are open dense in $X$,
and hence their intersection \begin{equation*}
X^{\circ}=\bigcap_{i=1}^lX_i=\{(x_1,\dots,x_l)\in X \mid x_i\ne 0
\text{ for any $i$}\} \end{equation*} is an open dense set in $X$,
and thus $X^{\circ}$ is irreducible. This implies that $T$ is
irreducible (and hence connected). Thus $T$ is a subtorus of
$T_l$. The assertion that $X$ is an equivariant affine embedding
of $T$ follows from $(1)$ and $(2)$.
\end{pf}
\begin{rem}
By Proposition \ref{above'}, we have that the ideal $I(X)$ is a
toric ideal in the sense of Definition \ref{tor.ideal}; more
precisely, $I(X)=I_{\mathcal{A}}$, where
${\mathcal{A}}=\{\rho_i,1\le i\le l\}$,  where
$\rho_i=\e_i\bigr|_T$ (here $K[T]$ is identified with $K[t_1^{\pm
1},\dots,t_d^{\pm 1}]$, and the character group $X(T)$ with
${\mathbb{Z}}^d$).
\end{rem}

\ni{\textbf{Summary:}} Summarizing we have established

(i) bijections between \{equivariant affine toroidal embeddings\},
\{toric ideals\} and \{binomial prime ideals\}.

(ii) a bijection between \{affine toric varieties with torus action
by $T$\}, \{finitely generated, saturated, sub semigroups of $M$
which generate $M$\}, and \{strongly convex rational polyhedral
cones inside $N_{\mathbb{R}}$\}.

\subsection{Orbit decomposition in affine toric varieties}\label{orbit} We
shall preserve the notation of \S \ref{comb}. As in \S \ref{comb},
let $X$ be an affine toric variety with a torus action by $T$. Let
$K[X]=K[S_\sigma]$. We shall denote $X$ also by $X_\sigma$.

 Let us first recall the definition of faces of
a convex polyhedral cone:
\begin{defn}
A face $\tau$ of $\sigma$ is a convex polyhedral sub cone of
$\sigma$ of the form $\tau=\sigma\cap u^{\perp}$ for some $u\in
\sigma^{\v}$, and is denoted $\tau<\sigma$.
\end{defn}

We have that $X_\tau$ is a principal open subset of $X_\sigma$,
namely, $$X_\tau=(X_\sigma)_u$$

\ni  Each face $\tau$ determines a (closed) point $P_\tau$ in
$X_\sigma$, namely, it is the point corresponding to the maximal
ideal in $K[X](=K[S_\sigma])$ given by the kernel of
$e_\tau:K[S_\sigma]\rightarrow K $, where for $u\in S_\sigma$, we
have
$$e_\tau(u)=\begin{cases}1,&{\mathrm{\ if\ }}u\in \tau^{\perp}\\
0,&{\mathrm{otherwise}} \end{cases}$$

\begin{rem}\label{coord} As a point in ${\mathbb{A}}^l$, $P_\tau$ may be identified with the
$l$-tuple with $1$ at the $i$-th place if $\chi_i$ is in
$\tau^{\perp}$, and $0$ otherwise (here, as in \S \ref{vanishing},
$\chi_i$ denotes the weight of the $T$-weight vector $y_i $ - the
class of $x_i$ in $K[X_\sigma]$).
\end{rem}

 Let $O_\tau$ denote
the $T$-orbit in $X_\sigma$ through $P_\tau$. In the sequel, we
shall refer to $P_\t$ as the \emph{center of} $O_\tau$. We have
the following orbit decomposition in $X_\sigma$:
$$\begin{gathered}
X_\sigma={\underset{\theta\le\sigma }{\cup}}\, O_\theta\\
{\overline{O_\tau}}={\underset{\theta\ge\tau }{\cup}}\, O_\theta\\
dim\,\tau + dim\,O_\tau=dim\,X_\sigma
\end{gathered}$$
(here, by dimension of a cone $\tau$, one means the vector space
dimension of the span of $\tau$).

See \cite{F}, \cite{toroidal} for details.

Thus $\tau\mapsto {\overline{O_\tau}}$ defines an order reversing
bijection between \{faces of $\sigma$\} and \{$T$-orbit closures
in $X_\sigma$\}. In particular, we have the following two extreme
cases:

\ni 1. If $\tau$ is the $0$-face, then $P_\tau=(1,\cdots,1)$ (as a
point in ${\mathbb{A}}^l$), and $O_\tau=T$, and is contained in
$X_\theta,\forall \theta<\sigma$. It is a dense open orbit.

\ni 2. If $\tau=\sigma$, then $P_\tau$ (as a point in
${\mathbb{A}}^l$) is the $l$-tuple $(0,\cdots,0)$ (we may suppose
that $\sigma$ spans $N_{\mathbb{R}}$); further,
$O_\tau=\{P_\tau\}$, and is the unique closed orbit.

For a face $\tau$, let us denote by $N_\tau$ the sublattice of $N$
 generated by the lattice points
of $\tau$. Let $N(\tau)=N/N_\tau$, and $M(\tau)$, the
${{\mathbb{Z}}}$-dual of $N(\tau)$.
For a face $\theta$ of $\sigma$ such that $\theta$ contains $\tau$
as a face, set
$$\theta_\tau:=(\theta+(N_\tau)_{{\mathbb{R}}})/(N_\tau)_{{\mathbb{R}}}$$
(here, for a lattice $L$, $L_{{\mathbb{R}}}=L\otimes
{\mathbb{R}}$). Then the collection
$\{\theta_\tau,\sigma\ge\theta\ge\tau\}$ forms the set of faces of
a cone in $ N(\tau)_{{\mathbb{R}}})$; we shall denote it by
$\sigma_\tau$.

\begin{lem}\label{closure} For a face $\tau<\sigma$, ${\overline{O_\tau}}$ gets
identified with the toric variety $Spec\,K[S_{\sigma_{\tau}}]$.
Further, $K[{\overline{O_\tau}}]=K[S_\sigma\cap\tau^{\perp} ]$.
\end{lem}
\begin{pf}
The first assertion follows from the description of the orbit
decomposition (and the definition of $\sigma_{\tau}$). For the
second assertion, we have
$$S_{\sigma_{\tau}}=\sigma_{\tau}^{\v}\cap
M(\tau)=\sigma^{\v}\cap(\tau^{\perp}\cap
M)=S_\sigma\cap\tau^{\perp}
$$(note that we have, $M(\tau)=\tau^{\perp}\cap M$).
\end{pf}

\section{Generalities on  finite distributive
lattices}\label{varlat} We shall now study a special class of
toric varieties, namely, the toric varieties associated to
distributive lattices. We shall first collect some definitions on
finite partially ordered sets. A partially ordered set is also
called a poset.
\begin{defn}
A finite poset $P$ is called\/ {\em bounded} if it has a unique
maximal, and a unique minimal element, denoted $\widehat{1}$ and
$\widehat{0}$ respectively.
\end{defn}
\begin{defn}
A totally ordered subset $C$ of a finite poset $P$ is called a\/
{\em chain}, and the number $\# C-1$ is called the\/ {\em length}
of the chain.
\end{defn}
\begin{defn}
A bounded poset $P$ is said to be\/ {\em graded} (or also
\emph{ranked}) if all maximal chains have the same length (note
that $\widehat{1}$ and $\widehat{0}$ belong to any maximal chain).
\end{defn}
\begin{defn}
Let $P$ be a graded poset. The length of a maximal chain in $P$ is
called the\/ {\em rank} of $P$.
\end{defn}
\begin{defn}\label{5.5}
Let $P$ be a graded poset. For $\l,\m\in P$ with $\l\ge\m$, the
graded poset $\{\t\in P\mid\m\le\t\le\l\}$ is called the\/ {\em
interval from $\m$ to $\l$}, and denoted by $[\m,\l]$. The rank of
$[\m,\l]$ is denoted by $l_{\m}(\l)$; if $\m=\widehat{0}$, then we
denote $l_{\m}(\l)$ by just $l(\l)$.
\end{defn}
\begin{rem}\label{level}
In the sequel, we refer to $l(\l)$ as \emph{the level of} $\l$ in
$\L$.
\end{rem}

\begin{defn}
Let $P$ be a graded poset, and $\l,\m\in P$, with $\l\ge\m$. The
ordered  pair $(\l,\m)$ is called a\/ {\em cover} (and we also say
that $\l$\/ {\em  covers} $\m$) if $l_{\m}(\l)=1$.
\end{defn}

\subsection{Generalities on distributive lattices}\label{s2}
\begin{defn} A {\em lattice\/} is a partially ordered set $(\mathcal{L},\le)$ such that,
for every pair of elements $x,y\in \cal{L}$, there exist elements
$x\v y$ and $x\w y$, called the\/ {\em join}, respectively the\/
{\em meet} of $x$ and $y$, defined by:
\begin{gather}
x\v y\ge x,\ x\v y\ge y,\text{ and if  }z\ge x \text{ and }
z\ge y,\text{ then } z\ge x\v y,\notag\\
x\w y\le x,\ x\w y\le y,\text{ and if  }z\le x \text{ and } z\le
y,\text{ then } z\le x\w y.\notag
\end{gather}
\end{defn}
It is easy to check that the operations $\v$ and $\w$ are
commutative and associative.
\begin{defn}
An element $z\in \cal{L}$ is called the\/ {\em zero} of $\cal{L}$,
denoted by\/ $\widehat{0}$, if $z\le x$ for all $x$ in $\cal{L}$.
An element $z\in \cal{L}$ is called the\/ {\em one} of $\cal{L}$,
denoted by\/ $\widehat{1}$, if $z\ge x$ for all $x$ in $\cal{L}$.
\end{defn}
\begin{defn}
Given a lattice ${\mathcal{L}}$, a subset ${\mathcal{L}}'\subset
{\mathcal{L}}$ is called a\/ {\em sublattice} of $\cal{L}$ if
$x,y\in{\mathcal{L}}'$ implies $x\w y\in{\mathcal{L}}'$, $x\v
y\in{\mathcal{L}}'$.
\end{defn}
\begin{defn}
Two lattices ${\mathcal{L}}_1$ and ${\mathcal{L}}_2$ are\/ {\em
isomorphic} if there exists a bijection
$\varphi:{\mathcal{L}}_1\to {\mathcal{L}}_2$ such that, for all
$x,y\in{\mathcal{L}}_1$,
\begin{equation*}
\varphi(x\v y)=\varphi(x)\v \varphi(y) \text{   and   }
\varphi(x\w y)=\varphi(x)\w \varphi(y).
\end{equation*}
\end{defn}
\begin{defn} A lattice is called\/ {\em distributive} if the following
identities hold:
\begin{align}
x\w (y\v z)&=(x\w y)\v (x\w z)\\
x\v (y\w z)&=(x\v y)\w (x\v z).
\end{align}
\end{defn}

\vs.2cm\ni\textbf{Example:}  The lattice of all subsets of the set
$\{1,2,\dots,n\}$ is a distributive lattice, denoted by
${\mathcal{B}}(n)$, and called the {\em Boolean algebra of rank\/}
$n$.

\begin{defn}
An element $z$ of a lattice $\L$ is called\/ {\em
join-irreducible} (respectively\/ {\em meet-irreducible}) if
$z=x\v y$ (respectively $z=x\w y$) implies $z=x$ or $z=y$. The set
of join-irreducible (respectively meet-irreducible) elements of
$\L$ is denoted by $J_{\L}$ (respectively $M_{\L}$), or just by
$J$ (respectively $M$) if
 no confusion is possible.
\end{defn}

\begin{defn}
The set $J_{\L}\cap M_{\L}$ of join and meet-irreducible elements
is denoted by $JM_{\L}$, or just $JM$ if no confusion is possible.
\end{defn}
\begin{defn}
A subset $I$ of a poset $P$ is called an {\em ideal} of $P$ if for
all $x,\,y\in P$,
$$x\in I\text{ and }y\le x\text{ imply }y\in I.$$
\end{defn}
\begin{thm}\label{5.10} (Birkhoff)
Let $\L$ be a distributive lattice with ${\hat{0}}$, and $P$ the
poset of its nonzero join-irreducible elements. Then $\L$ is
isomorphic to the lattice of finite ideals of $P$, by means of the
lattice isomorphism
$$\a\mapsto I_{\a}=\{\t\in P\mid \t\le\a\},\qquad \a\in\L.$$
\end{thm}
\begin{defn}
A quadruple of the form $(\t,\f,\t\v\f,\t\w\f)$, with $\t,\f\in\L$
non-comparable is called a\/ {\em diamond}, and is denoted by
$D(\t,\f,\t\v\f,\t\w\f)$ or also just $D(\t,\f)$. The pair
$(\t,\f)$ (respectively $(\t\v\f,\t\w\f)$)) is called the
\emph{skew} (respectively \emph{main}) diagonal of the diamond
$D(\t,\f)$.
\end{defn}
Denote
$$Q({\mathcal{L}})=\{(\tau,\varphi),\tau,\varphi{\mathrm{ noncomparable}}\}$$
In the sequel, $Q({\mathcal{L}})$ will also be denoted by just
$Q$.

The following Lemma is easily checked.
\begin{lem}
With the notations as above, we have

$(a)$ $J=\{\t\in\L\mid\text{there exists at most one cover of the
form } (\t,\l)\}$.

$(b)$ $M=\{\t\in\L\mid\text{there exists at most one cover of the
form } (\l,\t)\}$.
\end{lem}

For $\a\in\L$, let $I_\a$ be the ideal corresponding to $\a$ under
the isomorphism in Theorem \ref{5.10}.
\begin{lem}\label{cover}
Let $(\t,\l)$ be a cover in ${\mathcal{L}}$. Then $I_{\t}$ equals
$I_{\l}\dot\cup\{\beta\}$ for some $\beta\in J_{{\mathcal{L}}}$.
\end{lem}
\begin{pf}
If $\t\in J_{{\mathcal{L}}}$, then $\l$ is the unique element
covered by $\t$; it is clear that in this case that
$I_{\t}=I_{\l}\cup \{\t\}$.

Let then $\t\not\in J_{{\mathcal{L}}}$. Let $\l'$ be another
element covered by $\t$. Let $\f=\l\w\l'$.

\ni\textbf{Claim:} $(\l,\f),(\l',\f)$ are both covers.

We shall prove that $(\l,\f)$ is a cover, the proof being similar
for $(\l',\f)$. If possible, let us assume that there exists a
$\f'$ such that $\l>\f'>\f$. Then we have $\l'\le \f'\v\l'\le\t$.
The fact that $(\t,\l')$ is a cover implies that
$$\f'\v\l'=\t\leqno{(*)}$$ (note that $\f'\v\l'\ne \l'$; for, $\f'\v\l'=
\l'$ would imply that $\l'>\f'$, and this in turn would imply that
$\f(=\l\w\l')\ge\f'$ which would contradict the assumption that
$\f'>\f$). Also, we have,
$$\f=\l\w\l'\ge\f'\w\l'$$ (since $\l>\f'$), and
$$\f'\w\l'\ge \f$$ (since, both $\f'$ and $\l'$ are greater than
$\f$). Hence we obtain $$\f'\w\l'= \f\leqno{(**)}$$ Now in view of
(*) and (**), we obtain that for the pair $(\l,\f')$ (with
$\l>\f'$), we have $$\l\v\l'=\f'\v\l',\ \l\w\l'=\f'\w\l'$$ But
this is not possible, since ${\L}$ is a distributive lattice.
Hence our assumption is wrong and Claim follows.

 The above Claim together with induction on $l(\t)$
(the level of $\t$ (cf. Remark \ref{level})) implies that
$I_{\l}=I_{\f}\cup \{\beta\}$, for some $\beta\in
J_{{\mathcal{L}}}$. We have, $\beta\not\in I_{\l'}$ (note that
$\beta\in I_{\l'}$ would imply $I_{\l}\subseteq I_{\l'}$ (since $
\f\le\l'$) which is not true). We have that $I_{\l'}\cup
\{\beta\}$ is an ideal; for, any $\gamma<\beta$ (and $\gamma\neq
\beta$) is in fact $\le \phi$ (since $I_{\f}\cup \{\beta\}=I_{\l}$
is an ideal), and hence is in $I_{\l'}$. Further, $I_{\l'}\cup
\{\beta\}$ is clearly contained in $I_{\t}$ . Hence if $\theta$ is
the element of ${\mathcal{L}}$ corresponding to the ideal
$I_{\l'}\cup \{\beta\}$ (under the bijection given by Theorem
\ref{5.10}), then we have, $\l'< \theta\le \t$. Hence, we obtain
that $\theta=\t$ (since $(\t,\l')$ is a cover).

\vs.2cm\ni\textbf{Starting point of induction:} Let $\t$ be an
element of least length among all $\theta$'s such that $\theta$
covers an element $\theta'$. Then the above reasoning (especially,
the Claim) implies clearly that $\t\in J_{{\mathcal{L}}}$, and
therefore, $I_{\t}=I_{\l}\cup \{\t\}$.
\end{pf}
Combining the above Lemma with the fact that $I_{{\hat{1}}}$
equals $J_{{\mathcal{L}}}$, we have
\begin{cor}\label{chain}
Any maximal chain in (the ranked poset) ${\L}$ has cardinality
equal to $\# J_{{\mathcal{L}}}$.
\end{cor}
\section{The variety $X_{{\L}}$}\label{5.13} Consider the polynomial
algebra $K[X_\alpha,\alpha\in{\L}]$; let $I_{{\L}}$ be the ideal
generated by $\{X_\alpha X_\beta-X_{\alpha\v \beta}X_{\alpha\w
\beta}, \alpha,\beta\in{\L}\}$. Then one knows (cf.\cite{Hi}) that
$K[X_\alpha,\alpha\in{\L}]\,/I_{{\L}}$ is a normal domain; in
particular, we have that $I_{{\L}}$ is a  prime ideal. Let
$X_{{\L}}$ be the affine variety of the zeroes in $K^l$ of
$I_{{\L}}$ (here, $l=\# {\L}$).
 Then $X_{{\L}}$ is a affine normal variety defined by binomials. Let $T$ be as in Proposition \ref{binotoric}.
 Then by that Proposition, we have
 \begin{thm}\label{nortor}
 $X_{{\L}}$ is
 a (normal) toric variety for the action by $T$.
 \end{thm}

 We shall call $X_{{\L}}$ a \emph{Hibi toric variety}.
We shall now show that $\dim X_{{\L}}=\#J_{{\L}}$.
  We follow the
notation in \S \ref{varbino}.

Let ${\mathcal{I}}=\{(\t,\f,\t\v\f,\t\w\f)\mid (\t,\f)\in Q\}$,
where
\begin{equation*}
Q=\{(\t,\f)\mid \t,\f\in\L\text{ non-comparable}\}.
\end{equation*}

Let $T_l=(K^*)^l$, $\pi:X(T_l)\to X(T)$ be the canonical map,
given by restriction, and for $\h\in X(T_l)$, denote $\pi (\h)$ by
$\overline{\h}$. Let us fix a ${\mathbb{Z}}$-basis
$\{\h_\t\mid\t\in\L\}$ for $X(T_l)$. For a diamond
$D=(\t,\f,\t\v\f,\t\w\f)\in\cal{I}$, let
$\h_D=\h_{\t\v\f}+\h_{\t\w\f}-\h_{\t}-\h_{\f}$.

\begin{lem}\label{prev}
We have

$(1)$ $X(T)\simeq X(T_l)/\ker\pi$.

$(2)$ $\ker\pi$ is generated by $\{\h_D\mid D\in\cal{I}\}$.
\end{lem}
\begin{pf}
The restriction map  $\pi:X(T_l)\rightarrow X(T)$ is, in fact,
surjective, since $T$ is a subtorus of $T_l$. Now $(1)$ follows
from this. The assertion $(2)$ follows from the definition of $T$
(cf. Proposition \ref{binotoric}).
\end{pf}

Let $X(J_{{\mathcal{L}}})$ be the ${\mathbb{Z}}$-span of
$\{\chi_{\theta}, \theta\in J_{{\mathcal{L}}}\}$. As an immediate
consequence of the above Lemma, we have

\begin{cor}\label{isom}
$\pi$ maps $X(J_{{\mathcal{L}}})$ isomorphically onto its image.
\end{cor}

\begin{pf} The result follows since there does not exist a diamond contained
completely in $J_{{\mathcal{L}}}$.
\end{pf}


\begin{lem}\label{before'}
Let $\alpha\in\mathcal{L}$. Then $\overline{\chi}_\alpha$ is in
the $\mathbb{Z}$-span of $\{{\overline{\chi_\theta}},\theta\le
\alpha,\theta\in J_{\mathcal{L}}\}$. In particular,
$\{{\overline{\chi_\theta}},\,\theta\in J_{{\L}}\}$ generates
$X(T)$ as a $\mathbb{Z}$-module.
\end{lem}
\begin{pf} We shall prove the result by induction on $l(\alpha)$,
level of $\alpha$ (cf. Remark \ref{level}). If $l(\alpha)=0$, then
$\alpha={\hat{0}}$, and the result is clear. Let then
$l(\alpha)\ge 1$; this implies in particular that there exist
elements in ${\L}$ covered by $\alpha$. The result is clear if
$\alpha\in J_{{\mathcal{L}}}$.

Let then $\alpha\not\in J_{{\mathcal{L}}}$. This implies that
there exist $\alpha_1,\alpha_2$ both of which are covered by
$\alpha$. Then $\alpha_1,\alpha_2$ are non-comparable (clearly).
We have $\alpha_1\v\alpha_2=\alpha$; let
$\beta=\alpha_1\w\alpha_2$. In view of Lemma \ref{prev}(2), we
have,

 $${\overline{\chi}}_{\alpha}=
 {\overline{\chi}}_{\alpha_{1}}+{\overline{\chi}}_{\alpha_{2}}-{\overline{\chi}}_{\beta}$$
By induction, each term on the R.H.S. is in the $\mathbb{Z}$-span
of $\{{\overline{\chi_\theta}},\theta\le \alpha,\theta\in
J_{\mathcal{L}}\}$, and the result follows.
\end{pf}

Combining the above Lemma with Corollary \ref{isom}, we obtain
\begin{prop}
The set $\{{\overline{\h}}_\t\mid \t\in J_{\mathcal{L}}\}$ is a
${\mathbb{Z}}$-basis for $X(T)$.
\end{prop}

 Now, since $\dim X_{{\L}}=\dim T$, we obtain
\begin{thm}\label{dim}
The dimension of $X_{{\L}}$ is equal to  $\#J_{\L}$.
\end{thm}
Combining the above theorem with Lemma \ref{cover} and Corollary
\ref{chain}, we obtain
\begin{cor}\label{chain'}
\begin{enumerate}\item The dimension of $X_{{\L}}$ is equal to the cardinality of the set of elements in a
maximal chain in ${\L}$. \item Fix any chain
$\beta_1<\cdots<\beta_d$, $d$ being $\# J_{\L}$. Let
$\gamma_{i+1}$ be the element of $J_{\L}$ corresponding to the
cover $(\beta_{i+1},\beta_i)$, $i\ge 1$; set $\gamma_1=\beta_1$.
Then $J_{\L}=\{\gamma_1,\cdots,\gamma_d\}$
\end{enumerate}
\end{cor}
\begin{defn}
 For a finite distributive lattice $\L$, we call the
cardinality of $J_{\L}$  the {\em dimension} of $\L$, and we
denote it by $\dim\L$. If $\L'$ is a sublattice of $\L$, then the
{\em codimension} of $\L'$ in $\L$ is defined as $\dim\L-\dim\L'$.
\end{defn}

\begin{defn}\label{5.21}
A sublattice  $\L '$ of $\L$ is called an {\em embedded sublattice
of} $\L$ if
$$\t,\,\f\in\L,\quad\t\v\f,\,\t\w\f\in\L '\quad\Rightarrow\quad\t,\,\f\in\L '.$$
\end{defn}
Given a sublattice $\L '$ of $\L$, let us consider the variety
$X_{{\L}'}$, and consider the canonical embedding
$X_{{\L}'}\hookrightarrow{\mathbb{A}}(\L
')\hookrightarrow{\mathbb{A}}(\L)$ (here ${\mathbb{A}}(\L
')={\mathbb{A}}^{\#\L '}$,
${\mathbb{A}}(\L)={\mathbb{A}}^{\#\L}$).

\begin{prop}\label{5.22}
$X_{{\L}'}$ is a subvariety of $X_{{\L}}$ if and only if $\L '$ is
an embedded sublattice of $\L$.
\end{prop}
\begin{pf}
Under the embedding $X_{{\L}'}\hookrightarrow{\mathbb{A}}(\L)$,
$X_{{\L}'}$ can be identified with
$$\{(x_{\t})_{{\t\in\L}}\in{\mathbb{A}}{(\L)}\mid x_{\t}=0\text{ if }\t\not\in{\L} ',\ x_{\t}
x_{\f} =x_{\t\v\f}x_{\t\w\f}\text{ for }\t,\f\in{\L} '\text{
noncomparable }\}.$$ Let $\eta '$ be the center of the open dense
orbit in $X_{{\L}'}$ (cf. Remark \ref{coord},(1)); note that
$\eta'_{\t}\ne 0$ if and only if $\t\in{\L'}$. We have that
$X_{{\L}'}\subset X_{{\L}}$ if and only if $\eta '\in X_{{\L}}$.

Assume that $\eta '\in X_{{\L}}$. Let $\t$, $\f$ be two
noncomparable elements of $\L$ such that $\t\v\f$, $\t\v\f$ are
both in $\L '$. We have to show that $\t,\,\f\in\L '$. If
possible, let $\t\not\in\L '$. This implies $\eta_\t '=0$. Hence
either $\eta_{\t\v\f} '=0$, or $\eta_{\t\w\f} '=0$, since $\eta
'\in X_{{\L}}$. But this is not possible (note that $\t\v\f$,
$\t\w\f$ are in $\L '$, and hence $\eta_{\t\v\f} '$ and
$\eta_{\t\w\f} '$ are both nonzero).

Assume now that $\L '$ is an embedded sublattice. We have to show
that $\eta '\in X_{{\L}}$. Let $\t$, $\f$ be two non-comparable
elements of $\L$. The fact that $\L '$ is an embedded sublattice
implies that $\{\t,\f\}\subset {\L'}$ if and only if
$\{\t\v\f,\,\t\w\f\} \subset {\L'}$; further, when
$\t,\,\f,\,\t\v\f,\,\t\w\f$ are in ${\L'}$, we have
$$\eta_\t '\eta_\f '=\eta_{\t\v\f} '\eta_{\t\w\f} '$$
Thus $\eta '$ satisfies the defining equations of $X_{{\L}}$, and
hence $\eta '\in X_{{\L}}$.
\end{pf}



\section{Cone and dual cone of $X_{{\mathcal{L}}}$:}\label{cone} In this section we shall determine the cone
 and the dual cone of $X_{{\mathcal{L}}}$, ${\mathcal{L}}$ being a finite distributive lattice.
  As in the previous sections, denote the poset of join-irreducibles by
$J_{{\mathcal{L}}}$; let $d=\#\,J_{{\mathcal{L}}}$. Denote
$J_{{\mathcal{L}}}=\{\beta_{1},\cdots,\beta_{d}\}$. Let $T$ be a
$d$-dimensional torus. Identifying $T$ with $(K^*)^d$, let
$\{f_\beta,\beta\in
  J_{{\mathcal{L}}}\}$ denote the standard ${\mathbb{Z}}$-basis  for
$X(T)$, namely, for $u\in T, u=(u_\beta, \beta \in
J_{{\mathcal{L}}})$, $f_\beta(u)=u_\beta$.  Then $K[X(T)]$ may be
identified with $K[u_{\beta_{1}}^{\pm 1},\cdots,u_{\beta_{d}}^{\pm
1}]$, the ring of Laurent polynomials. If $\chi\in X(T)$, say
$\chi=\sum\,a_\beta f_\beta$, then as an element of $K[X(T)]$,
$\chi$ will also be denoted by $\prod\,u_\beta^{a_\beta}$.

Denote by ${\mathcal{I}}(J_{{\mathcal{L}}})$ the poset of ideals
of $J_{{\mathcal{L}}}$. For $A\in
{\mathcal{I}}(J_{{\mathcal{L}}})$, set
$$f_A:={\underset{z\in A}{\sum}} \,f_z$$
Recall (cf. Theorem \ref{5.10}) the bijection ${\L}\rightarrow
{\mathcal{I}}(J_{{\mathcal{L}}}),\theta\mapsto A_\theta:=\{\tau\in
 J_{{\mathcal{L}}}\,|\,\tau\le\theta \}$.
 \begin{lem}
$\{f_{A_{\theta}}, \theta\in J_{{\mathcal{L}}}\}$ is a
$\mathbb{Z}$-basis for $X(T)$.
 \end{lem}
 \begin{pf}
Take a total order on $J_{{\mathcal{L}}}$ extending the partial
order. Then the matrix expressing $\{f_{A_{\theta}}, \theta\in
J_{{\mathcal{L}}}\}$ in terms of the basis $\{f_{\theta},
\theta\in J_{{\mathcal{L}}}\}$ is easily seen to be triangular
with diagonal entries equal to 1. The result follows from this.
 \end{pf}

 Let ${\mathcal{A}}=\{f_A,A\in{\mathcal{I}}(J_{{\mathcal{L}}})\}$.
 Then as a consequence of the above Lemma, we obtain
 \begin{cor}\label{above}
$\mathcal{A}$ generates $X(T)$; in particular,
$d_{{\mathcal{A}}}=d,\ d_{{\mathcal{A}}}$ being as in Proposition
\ref{zar}.
 \end{cor}

  For  $A\in {\mathcal{I}}(J_{{\mathcal{L}}})$,
denote by ${\mathbf{m}}_A$ the monomial:
$${\mathbf{m}}_A:={\underset{\tau\in A}{\prod}} \,u_\tau$$ in  $K[X(T)]$. If
$\alpha$ is the element of ${\L}$ such that $I_\alpha=A$ (cf.
Theorem \ref{5.10}), then we shall denote ${\mathbf{m}}_A$ also by
${\mathbf{m}}_\alpha$. Consider the surjective algebra map
$$F:K[X_\alpha,\alpha\in{\mathcal{L}}]\rightarrow
K[{\mathbf{m}}_A,A\in{\mathcal{I}}(J_{{\mathcal{L}}})]\,(\subset
K[X(T)]),\,X_\alpha\mapsto {\mathbf{m}}_A,\ A=I_\alpha
$$ Then with notation as in \S \ref{toric} and Definition \ref{tor.ideal}, we have, $F=\widehat{\pi_{{\mathcal{A}}}},ker\,F=I_{\mathcal{A}}$.
Further, Proposition \ref{zar} and Corollary \ref{above} imply the
following
\begin{prop}\label{zar'}
${\mathcal{V}}(I_{{\mathcal{A}}})=Spec\,K[{\mathbf{m}}_A,A\in{\mathcal{I}}(J_{{\mathcal{L}}})]$
and is of dimension $d (=\#J_{{\mathcal{L}}})$ (here,
${\mathcal{V}}(I_{\mathcal{A}})$ is as in Proposition \ref{zar})
\end{prop}

Let us denote ${\mathcal{V}}(I_{\mathcal{A}})$ by $Y$.
\begin{lem}  Kernel of $F$ is generated by $\{X_\alpha X_\beta-X_{\alpha\v \beta}X_{\alpha\w
\beta}, \alpha,\beta\in\L\}$
\end{lem}
\begin{pf}
We have (in view of Lemma \ref{cover}) that if $(\alpha,\l)$ is a
cover, then $I_\alpha$ equals $I_\l\dot\cup\{\beta\}$ for some
$\beta\in J_{\L}$. A repeated application of this result implies
that if $\gamma\le\alpha$, and $m=l(\alpha)-l(\gamma)$ (here,
$l(\beta)$ denotes the level of $\beta$ (cf. Remark \ref{level})),
then there exist $\alpha_1\cdots,\alpha_m$ in $J_{{\L}}$ such that
$I_\alpha\,\setminus\,I_\gamma=\{\alpha_1\cdots,\alpha_m\}$. Let
now $\beta,\beta'$ be two non-comparable elements in $\L$. Let
$\alpha=\beta\v\beta',\gamma=\beta\w\beta'$. Let
$l(\beta)-l(\gamma)=r,l(\beta')-l(\gamma)=s$. Then there exist
$\beta_1,\cdots,\beta_r$, and $\beta'_1,\cdots,\beta'_s$ in
$J_{{\L}}$ such that
$$\begin{gathered} I_\beta\,\setminus\,I_\gamma=\{\beta_1,\cdots,\beta_r\}\\
I_{\beta'}\,\setminus\,I_\gamma=\{\beta'_1,\cdots,\beta'_s\}\\
I_\alpha\,\setminus\,I_\beta=\{\beta'_1,\cdots,\beta'_s\}\\
I_\alpha\,\setminus\,I_{\beta'}=\{\beta_1,\cdots,\beta_r\}\end
{gathered}$$ Hence we obtain
$$\begin{gathered}{\mathbf{m}}_\beta={\mathbf{m}}_\gamma
u_{\beta_{1}}\cdots u_{\beta_{r}}\\
{\mathbf{m}}_{\beta'}={\mathbf{m}}_\gamma u_{\beta'_{1}}\cdots
u_{\beta'_{s}}\\
{\mathbf{m}}_{\alpha}={\mathbf{m}}_\beta u_{\beta'_{1}}\cdots
u_{\beta'_{s}}\\
{\mathbf{m}}_{\alpha}={\mathbf{m}}_{\beta'}u_{\beta_{1}}\cdots
u_{\beta_{r}}
\end{gathered}$$ From this it follows that $${\mathbf{m}}_{\alpha}{\mathbf{m}}_{\gamma}={\mathbf{m}}_{\beta}
{\mathbf{m}}_{\beta'}$$ Thus for each diamond in $\L$, i.e., a
quadruple $(\beta,\beta',\beta\v\beta',\beta\w\beta')$, we have
$X_\beta X_{\beta'}-X_{\beta\v\beta'}X_{\beta\w\beta'}$ is in the
kernel of the surjective map $F$. Hence $F$ factors through
$K[X_{\L}]$; hence, we obtain  closed immersions (of irreducible
varieties):
$$Y\hookrightarrow X_{\L}\hookrightarrow {\mathbb{A}}^{\#{\L}}$$ ($Y$ being
${\mathcal{V}}(I_{\mathcal{A}})$). But dimension considerations
 imply that
$Y=X_{\L}$ (note that in view of Proposition \ref{zar'},
$dim\,Y=d=dim\,X_{\L}$ (cf. Theorem \ref{dim})), and the result
follows.
\end{pf}
As an immediate consequence, we obtain (as seen in the proof of
the above Lemma)
\begin{thm}\label{main1} We have an isomorphism
$K[X_{\L}]\cong
K[{\mathbf{m}}_A,A\in{\mathcal{I}}(J_{{\mathcal{L}}})]$.
\end{thm}

 Denote $M:=X(T)$, the character group pf
$T$. Let $N={\mathbb{Z}}$-dual of $M$, and $\{e_y,y\in
J_{{\mathcal{L}}}\}$ be the basis of $N$ dual to $\{f_z,z\in
J_{{\mathcal{L}}}\}$. As above,  for $A\in
{\mathcal{I}}(J_{{\mathcal{L}}})$, let
$$f_A:={\underset{z\in A}{\sum}} \,f_z$$ Let
$V=N_{{\mathbb{R}}}(=N\otimes_{{\mathbb{Z}}} {\mathbb{R}})$. Let
$\sigma \subset V$ be the cone such that
$X_{{\mathcal{L}}}=X_\sigma$. Let $\sigma^{\v}\subset V^*$ be the
cone dual to $\sigma$. Let $S_\sigma=\sigma^{\v}\cap M$, so that
$K[X_{{\mathcal{L}}}]$ equals the semi group algebra
$K[S_\sigma]$.

As an immediate consequence of Theorem \ref{main1}, we have
\begin{prop}\label{semi}
The semigroup $S_\sigma$ is generated by $\{f_A, A\in
{\mathcal{I}}(J_{{\mathcal{L}}})\}$.
\end{prop}

Let $M(J_{{\mathcal{L}}})$ be the set of maximal elements in the
poset $J_{{\mathcal{L}}}$. Let $Z(J_{{\mathcal{L}}})$ denote the
set of all covers in the poset $J_{{\mathcal{L}}}$ (i.e., $(z,z'),
z>z'$ in the poset $J_{{\mathcal{L}}}$, and there is no other
element $y\in J_{{\mathcal{L}}}$ such that $z>y>z'$). For a cover
$(y,y')\in Z(J_{{\mathcal{L}}})$, denote
$$v_{y,y'}:=e_{y'}-e_{y}$$

\begin{prop} The cone $\sigma$ is
generated by $\{e_{z},z\in M(J_{{\mathcal{L}}}),\,v_{y,y'},
(y,y')\in Z(J_{{\mathcal{L}}})\}$.
\end{prop}
\begin{pf}
Let $\theta$ be the cone generated by $\{e_{z},z\in
M(J_{{\mathcal{L}}}),\,v_{y,y'}, (y,y')\in
Z(J_{{\mathcal{L}}})\}$. Then clearly any $u$ in $S_\sigma$ is
non-negative on the generators of $\theta$, and hence
$\sigma^{\v}\subseteq \theta^{\v}$. We shall now show that
$\theta^{\v}\subseteq\sigma^{\v}$, equivalently, we shall show
that $S_\theta\subseteq S_\sigma$. Let $f\in M$, say,
$f={\underset{z\in J_{{\mathcal{L}}}}{\sum}}\,a_zf_z$. Then it is
clear that $f$ is in $S_\theta$ if and only if $$a_z\ge
0,{\mathrm{\ for\ }}z{\mathrm{\  maximal\ in\ }}
J_{{\mathcal{L}}},{\mathrm{\ and\ }} a_x\ge a_y, {\mathrm{\ for\
}}x,y\in J_{{\mathcal{L}}}, x<y\leqno{(*)}$$

\ni\textbf{Claim:} Let $f\in M$. Then $f$ has property (*) if and
only if $f={\underset{A\in
{\mathcal{I}}(J_{{\mathcal{L}}})}{\sum}}c_Af_A$,

\ni $c_A\in {\mathbb{Z}}_+$.

Note that Claim implies that $S_\theta\subseteq S_\sigma$, and the
required result follows.

\ni\textbf{Proof of Claim:} The implication $\Leftarrow$ is clear.

\ni\textbf{The implication $\Rightarrow$:} Let $f\in S_\theta$,
say, $f={\underset{z\in J_{{\mathcal{L}}}}{\sum}}\,a_zf_z$. The
hypothesis that $f$ has property (*) implies that $a_z$ being
non-negative for $z$ maximal in $J_{{\mathcal{L}}}$, $a_x$ is
non-negative for all $x\in J_{{\mathcal{L}}}$. Thus
$f={\underset{z\in J_{{\mathcal{L}}}}{\sum}}\,a_zf_z,a_z\in
{\mathbb{Z}}_+,z\in J_{{\mathcal{L}}}$. Further, the property in
(*) that $a_x\ge a_y$ if $x<y$ implies that $$\{x\,|\,a_x\ne 0\}$$
is an ideal in $J_{{\mathcal{L}}}$. Call it $A$; in the sequel we
shall denote $A$ also by $I_f$. Let $m=\,$min\{$a_x,x\in A$\}.
Then either $f=mf_A$ in which case the Claim follows, or
$f=mf_A+f_1$, where $f_1$ is in $S_\theta$ (note that $f_1$ also
has property (*)); further, $I_{f_{1}}$ is a proper subset of $A$.
Thus proceeding we arrive at positive integers $m_1,\cdots,m_r$,
elements $f_1,\cdots,f_r$ (in $S_\theta$) with $A_i:=I_{f_{i}}$,
an ideal in $J_{{\mathcal{L}}}$, and a proper subset of
$A_{i-1}:=I_{f_{i-1}}$ such that
$$f=mf_A+m_1f_{A_{1}}+\cdots+m_rf_{A_{r}}$$ (In fact, at the last step, we have,
$f_r=m_rf_{A_{r}}$.) The Claim and hence the required result now
follows.
\end{pf}
\subsection{Analysis of faces of $\sigma$}\label{anal} We shall concern ourselves
 just with the closed points in $X_{\L}$.
So in the sequel, by a point in $X_{\L}$, we shall mean a closed
point. Let $\tau$ be a face of $\sigma$. Let $P_{\tau}$ be the
distinguished point of $O_{\tau}$ with the associated maximal
ideal being the kernel of the map
$$\begin{gathered}K[S_\sigma]\rightarrow K,\\
u\in S_\sigma,u\mapsto\begin{cases} 1,&\mathrm{\ if\ }
u\in\tau^{\perp}\\
0,&\mathrm{\ otherwise}
\end{cases}
\end{gathered}$$

Then for a point $P\in X_{\L}$ (identified with a point in
${\mathbb{A}}^l$), denoting by $P(\alpha)$, the $\alpha$-th
co-ordinate of $P$, we have, $$P_{\tau}(\alpha)=\begin{cases}
1,&\mathrm{\ if\ }
f_{I_{\alpha}}\in\tau^{\perp}\\
0,&\mathrm{\ otherwise}
\end{cases}$$
With notation as above, let
$$D_{\tau}=\{\alpha\in{\L}\,|\,P_{\tau}(\alpha)\neq 0\}$$ We have,
\begin{lem}\label{embed}
$D_{\tau}$ is an embedded sublattice of ${\L}$.
\end{lem}
\begin{pf}
Let $\theta,\delta$ be a pair of non-comparable elements in
$D_{\tau}$. The fact that $P_{\tau}(\theta),P_{\tau}(\delta)$ are
non-zero, together with the diamond relation $x_{\theta}
x_{\delta}=x_{\theta\v\delta} x_{\theta\w\delta}$ implies that

\ni $P_{\tau}(\theta\v\delta),P_{\tau}(\theta\w\delta)$ are again
non-zero, and are equal to $1$ (note that any non-zero co-ordinate
in $P_{\tau}$ is in fact equal to $1$). Thus, $D_{\tau}$ is a
 sublattice of ${\L}$.

The above reasoning implies that if $\theta,\delta$ in $D_{\tau}$
form the main diagonal in a diamond $D$ in ${\L}$, then
$P_{\tau}(\alpha),P_{\tau}(\beta)$ are non-zero, and are equal to
$1$, where $\alpha,\beta$ form the skew diagonal of $D$. Hence
$\alpha,\beta$ are in $D$, and hence $D\subset D_\tau$. Thus we
obtain that $D_{\tau}$ is an embedded sublattice of ${\L}$.
\end{pf}

Conversely, we have
\begin{lem}
Let ${\mathcal{D}}$ be an embedded sublattice in ${\L}$. Then
${\mathcal{D}}$ determines a unique face $\tau$ of $\sigma$ such
that $D_{\tau}$ equals ${\mathcal{D}}$.
\end{lem}
\begin{pf}
Denote $P$ to be the point in ${\mathbb{A}}^{\#{\L}}$ with
$$P(\alpha)=\begin{cases} 1,&{\mathrm{\ if\ }}
{\alpha}\in {\mathcal{D}}\\
0,&\mathrm{\ otherwise}
\end{cases}$$ Then $P$ is in ${\L}$ (since ${\mathcal{D}}$ is an embedded sublattice in ${\L}$,
the co-ordinates of $P$ satisfy all the diamond relations). Set
$$u={\underset{\alpha\in
{\mathcal{D}}}{\sum}}\,f_{I_{\alpha}},\,\tau=\sigma\cap
u^{\perp}$$ Then clearly, $P=P_{\tau}$ and
$D_{\tau}={\mathcal{D}}$
\end{pf}
Thus in view of the two Lemmas above, we have a bijection
$$\{\mathrm{\ faces\ of\
}\sigma\}\,{\buildrel{bij}\over{\leftrightarrow}}\, \{\mathrm{\
embedded\ sublattices \ of\ }{\L}\}$$
\begin{prop}\label{close}
Let $\tau$ be a face of $\sigma$. Then we have
${\overline{O_{\tau}}}=X_{D_{\tau}}$.
\end{prop}
\begin{pf}
Recall (cf. Lemma \ref{closure}) that
$K[{\overline{O_{\tau}}}]=K[S_\sigma\cap\tau^{\perp}]$. From the
description of $P_{\tau}$, we have that $\tau^{\perp}$ is the span
of $\{f_{I_{\alpha}},\alpha\in D_{\tau}\}$; hence
$$S_\sigma\cap\tau^{\perp}=\{{\underset{\alpha\in
D_{\tau}}{\sum}}\,c_{\alpha}f_{I_{\alpha}},c_{\alpha}\in\mathbb{Z}^+\}\
\leqno{(*)}$$ Thus we obtain
$$S_{\sigma}\cap\tau^{\perp}=S_{\sigma_{\tau}}\ \leqno{(**)}$$ (here, $\sigma_{\tau}$ is as in Lemma
\ref{closure}). On the other hand, if $\eta$ is the cone
associated to the toric variety $X_{D_{\tau}}$, then by
Proposition \ref{semi} we have  that $S_\eta$ is the semigroup
generated by $\{f_{I_{\alpha}},\alpha\in D_{\tau}\}$. Hence we
obtain that $\eta=\sigma_\tau$ (in view of (*) and (**)). The
required result now follows.
\end{pf}
\section{Generation of the cotangent space at $P_{\tau}$}\label{tan-cot} For $\alpha\in{{\L}}$,
let us denote the image of $X_\alpha$ in $R_{{\L}}$ (under the
surjective map

\ni $K[X_\alpha,\alpha\in{{\L}}]\rightarrow R_{{\L}}$) by
$x_\alpha$. Let $R_\tau=K[X(D_\tau)]$, the co-ordinate ring of
$X(D_\tau)$, and $\pi_\tau:R_{{\L}}\rightarrow R_\tau$ be the
canonical surjective map induced by the closed immersion
$X(D_\tau)\hookrightarrow X_{{\L}}$; clearly, kernel of $\pi_\tau$
is generated by $\{x_\theta,\theta\in {\L}\,\setminus\,D_\tau\}$.
Set
$$F_{\alpha}=\begin{cases} x_\alpha,&{\mathrm{\ if\ }}
{\alpha}\not\in D_{\tau}\\
1-x_\alpha,&{\mathrm{\ if\ }}\alpha\in D_{\tau}
\end{cases}$$ Denoting by $M_{\tau}$, the maximal ideal in
$R_{\L}$ corresponding to $P_{\tau}$, we have (cf. \S \ref{anal})
\begin{lem}
The ideal $M_{\tau}$ is generated by
$\{F_{\alpha},\alpha\in{\L}\}$.
\end{lem}
\subsubsection{A set of generators for the cotangent space $M_{\tau}/M^2_{\tau}$}
For $F\in M_\tau$, let ${\overline{F}}$ denote the class of $F$ in
$M_{\tau}/M^2_{\tau}$.
\begin{lem}\label{two}
Fix a maximal chain $\Gamma$ in $D_{\tau}$. For any $\beta\in
D_{\tau}\,\setminus\,\Gamma$, we have that in
$M_{\tau}/M^2_{\tau}$, ${\overline{F}}_\beta$ is in the span of
$\{{\overline{F}}_\gamma,\gamma\in\Gamma\}$
\end{lem}
\begin{pf}
Fix a $\beta\in D_{\tau}\,\setminus\,\Gamma$; denote by
$l_\tau(\beta)$, the level (cf. Remark \ref{level}) of $\beta$
considered as an element of the distributive lattice $D_{\tau}$.
We shall prove the result by induction on $l_\tau(\beta)$. If
$l_\tau(\beta)=0$, then $\beta$ coincides with the (unique)
minimal element of $D_{\tau}$, and there is nothing to prove. Let
then $l_\tau(\beta)\ge 1$. Let $\beta'$ in $D_\tau$ be such that
$\beta$ covers $\beta'$. Then
$I_\beta\,\setminus\,I_{\beta'}=\{\theta\}$, for an unique
$\theta\in J_{D_{\tau}}$ (cf. Lemma \ref{cover}). Then there
exists a unique cover $(\gamma,\gamma')$ in $\Gamma$ such that
$I_\gamma\,\setminus\,I_{\gamma'}=\{\theta\}$ (cf. Corollary
\ref{chain'},(2)).

\ni\textbf{Claim:} $F_\beta-F_{\beta'}\equiv
F_\gamma-F_{\gamma'}(mod\,M_{\tau})$.

First observe that the fact that $\theta$ belongs to
$I_\beta,I_\gamma$ and does not belong to $I_{\beta'},I_{\gamma'}$
implies $$\gamma'\not\ge\beta,\,\beta'\not\ge\gamma \leqno{(*)}$$
We now divide the proof of the Claim into the following two cases.

\ni {\textbf{Case 1:}} $\gamma'<\beta$.

This implies that $\gamma'$ is in fact less than $\beta'$, and
$\gamma<\beta$ (since,
$I_\beta\,\setminus\,I_{\beta'}=I_\gamma\,\setminus\,I_{\gamma'}=\{\theta\}$,
and $\theta\not\in I_{\gamma'}$); this in turn implies that
$\beta'\not\le\gamma$ (for, otherwise, $\beta'\le \gamma$ would
imply $\gamma'<\beta'<\gamma$, not possible, since,
$(\gamma,\gamma')$ is a cover). Let $D$ be the diamond having
$(\gamma,\beta')$ as the skew diagonal. Now $\gamma'$ being less
than both $\gamma$ and $\beta'$, we get that
$\gamma'\le\beta'\w\gamma$, and in fact equals $\beta'\w\gamma$
(note that the fact that $\gamma'\le\beta'\w\gamma<\gamma$
together with the fact that $(\gamma,\gamma')$ is a cover in
$D_{\tau}$ implies that $\gamma'=\beta'\w\gamma$). In a similar
way, we obtain that $\beta=\beta'\v\gamma$. Now the diamond
relation $x_{\beta'}x_{\gamma}=x_{\beta}x_{\gamma'}$ implies the
claim in this case (by definition of $F_\xi$'s).

\ni {\textbf{Case 2:}} $\gamma'\not\le\beta$.

This implies that $\gamma\not\le\beta$.

If $\gamma>\beta$, then clearly
$$\gamma=\gamma'\v\beta,\beta'=\gamma'\w\beta$$ Claim in this case
follows as in Case 1.

Let then $\gamma\not\ge\beta$. Then we have that
$(\gamma,\beta),(\gamma',\beta')$ are non-comparable pairs (note
that $\gamma'\le\beta'$ if and only if $\gamma'<\beta$). Denote
$\delta=\gamma\w\beta,\delta'=\gamma'\w\beta'$. The facts that
$$I_\delta=I_\gamma\cap I_\beta,I_{\delta'}=I_{\gamma'}\cap I_{\beta'},
I_\gamma=I_{\gamma'}\cup\{\theta\},I_\beta=I_{\beta'}\cup\{\theta\}$$
imply that $$I_\delta=I_{\delta'}\cup\{\theta\}$$ Hence we obtain
$$\begin{gathered} I_{\gamma'}\cap I_\delta=I_{\gamma'}\cap
I_{\delta'}=I_{\delta'};I_{\gamma'}\cup I_\delta=I_\gamma\\
I_{\beta'}\cap I_\delta=I_{\beta'}\cap
I_{\delta'}=I_{\delta'};I_{\beta'}\cup I_\delta=I_\beta
\end{gathered}$$ Thus we obtain that
$$\begin{gathered}
\gamma'\w\delta=\delta',\gamma'\v\delta=\gamma\\
\beta'\w\delta=\delta',\beta'\v\delta=\beta
\end{gathered}$$ Considering the diamonds with skew diagonals $(\gamma',\delta),
(\beta',\delta)$respectively, the diamond relations
$$x_{\gamma'}x_{\delta}=x_{\gamma}x_{\delta'},
x_{\beta'}x_{\delta}=x_{\beta}x_{\delta'}$$imply that in
$M_{\tau}/M^2_{\tau}$, we have the following relations
$${\overline{F}}_\gamma-{\overline{F}}_{\gamma'}={\overline{F}}_\delta-{\overline{F}}_{\delta'};
{\overline{F}}_\beta-{\overline{F}}_{\beta'}={\overline{F}}_\delta-{\overline{F}}_{\delta'}$$
Hence we obtain that $F_\beta-F_{\beta'}\equiv
F_\gamma-F_{\gamma'}(mod\,M_{\tau})$ as required.

This completes the proof of the Claim. Note that Claim implies the
required result (by induction on $l_\tau(\beta)$). (Note that when
$l_\tau(\beta)=1$, then $\beta'$ is the (unique) minimal element
of $D_{\tau}$, and hence $\beta'\in \Gamma$. The result in this
case follows from the Claim).
\end{pf}
Under the (surjective) map $\pi_\tau:R_{\L}\rightarrow R_\tau$,
denote $\pi_\tau(M_\tau)$ by $M'_\tau$; then $\pi_\tau$ induces a
surjection $\pi_\tau:M_\tau/M^{2}_\tau\rightarrow
M'_\tau/M{'}^{2}_\tau$.
\begin{cor}\label{dtau}
 $\{\pi_\tau({\overline{F}}_\gamma),\gamma\in\Gamma\}$ is a basis for
$M'_\tau/M{'}^{2}_\tau$.
\end{cor}
\begin{pf}
We have that $dim\,M'_\tau/M{'}^{2}_\tau\ge
dim\,X(D_\tau)=\#\Gamma$ (cf. Corollary \ref{chain'}(1)). On the
other hand, the above Lemma implies that
$dim\,M'_\tau/M{'}^{2}_\tau\le \#\Gamma$, and the result follows.
\end{pf}
\begin{lem}\label{one}
Let $\alpha\in{\L}\,\setminus\,D_{\tau}$ be such that there exists
an element $\beta\in D_{\tau}$ and a diamond $D$ in ${\L}$ such
that

(1) $(\alpha,\beta)$ is a diagonal (main or skew) in $D$.

(2) $D\cap D_{\tau}=\{\beta\}$.

\ni Then $F_\alpha\in M^2_{\tau}$.
\end{lem}
\begin{pf} Let us denote the remaining vertices of the diamond by
$\theta,\delta$. Writing the diamond relation $x_\alpha
x_\beta=x_\theta x_\delta$ in terms of the $F_\xi$'s, we have,
$$x_\alpha x_\beta-x_\theta x_\delta=F_\alpha (1-F_\beta)-F_\theta F_\delta$$
Hence in $R_{\L}$, we have, $$F_\alpha=F_\alpha F_\beta+F_\theta
F_\delta$$ The required result follows from this.
\end{pf}

\subsubsection{The set $E_{\tau}$}\label{sub} Define
$$E_{\tau}:=\{\alpha\in{\L},\alpha\mathrm{\ as\ in\ Lemma\
\ref{one}}\}$$In the sequel, we shall refer to an element
$\alpha\in E_{\tau}$ as an $E_{\tau}$-element.

\subsubsection{The equivalence relation:}\label{equ} For two distinct elements $\theta,\delta\in
{\L}\,\setminus\,D_\tau$, we say \emph{$\theta$ is equivalent to
$\delta$} (and denote it as $\theta\sim\delta$) if there exists a
sequence $\theta=\theta_1,\cdots,\theta_n=\delta$ in
${\L}\,\setminus\,D_\tau$  such that $(\theta_i,\theta_{i+1})$
forms one side of a diamond in ${\L}$ whose other side is
contained in $D_\tau$. For $\theta\in {\L}\,\setminus\,D_\tau$, we
shall denote by $[\theta]$, the set of all elements of
${\L}\,\setminus\,D_\tau$ equivalent to $\theta$, if there exists
such a diamond as above having $\theta$ as one vertex; if no such
diamond exists, then $[\theta]$ will denote
  the singleton set $\{\theta\}$. Clearly, for all $\theta$ in a given
equivalence class, $F_\theta$ (mod$\,M_\tau^2$) is the same (by
consideration of diamond relations), and we shall denote it by
${\overline{F}}_{[\theta]}$ or also just ${\overline{F}}_\theta$.
Note also that in view of Lemma \ref{one}, ${\overline{F}}_\theta
=0$ in $M_{\tau}/M^2_{\tau}$, if $[\theta]\cap E_\tau \ne\emptyset
$. We shall refer to $[\theta]$ as a \emph{$E_\tau$-class or a
non$\,E_\tau$-class} according as $[\theta]\cap E_\tau$ is
non-empty or empty.

\vs.2cm\ni\textbf{The sublattice $\Lambda_{\tau}(\Gamma)$:} Fix a
chain $\Gamma$ in $D_{\tau}$. Let $\Lambda_{\tau}(\Gamma)$  denote
the union of all the maximal chains in ${\L}$ containing $\Gamma$.
Note that $\alpha$ in ${\L}$ is in $\Lambda_{\tau}(\Gamma)$ if and
only if $\alpha$ is comparable to every $\gamma\in\Gamma$.
\begin{lem}\label{three}
$\Lambda_{\tau}(\Gamma)$ is a sublattice of ${\L}$.
\end{lem}
\begin{pf}
Let $(\theta,\delta)$ be a pair of non-comparable elements in
$\Lambda_{\tau}(\Gamma)$. Then for any $\gamma\in\Gamma$, we have
(by definition of $\Lambda_{\tau}(\Gamma)$) that either $\gamma$
is less than both $\theta$ and $\delta$ or greater than both
$\theta$ and $\delta$; in the former case, $\gamma$ is less than
both $\theta\v\delta$ and $\theta\w\delta$, and in the latter
case, $\gamma$ is greater than both $\theta\v\delta$ and
$\theta\w\delta$.
\end{pf}
\begin{rem}
$\Lambda_{\tau}(\Gamma)$ need not be an embedded sublattice:

Take ${\L}$ to be the distributive lattice consisting of $\{(i,j),
1\le i,j\le 4\}$ with the partial order $(a_1,a_2)\ge
(b_1,b_2)\Leftrightarrow a_r\ge b_r, r=1,2$. Take $$D_\t:=\{(i,j),
2\le i,j\le 3\},\ \Gamma:=\{(2,2),(3,2),(3,3)\}$$ Then
$$\Lambda_{\tau}(\Gamma)=\Gamma\cup I_1\cup I_2$$ where
$I_1=\{(1,1),(1,2),(2,1)\},I_2=\{(3,4),(4,3),(4,4)\}$. Consider
$x=(2,4),y=(4,2)$; then $x\v y=(4,4),x\w y=(2,2)$. Now $x\v y,x\w
y$ are in $\Lambda_{\tau}(\Gamma)$, but $x,y$ are not in
$\Lambda_{\tau}(\Gamma)$. Thus $\Lambda_{\tau}(\Gamma)$ is not an
embedded sublattice.
\end{rem}
\begin{prop}\label{four} Let
$\theta_0\in {\L}\,\setminus\, D_\tau$. Then $\theta_0\sim\mu$,
for some $\mu\in \Lambda_{\tau}(\Gamma)\cup E_\tau$.
\end{prop}
\begin{pf}
  Let us denote ${\L}'_\tau={\L}\,\setminus\, D_\tau$.
  By induction, let us suppose that the result holds for all $\theta\in{\L}'_\tau,\theta>\theta_0$;
  we shall see that the proof for the case when $\theta_0$ is a maximal
element in ${\L}'_\tau$  (the starting point of induction) is
included in the proof for a general $\theta_0$.

  If $\theta_0\in\Lambda_{\tau}(\Gamma)$, then there is nothing to prove.
  Let then
$\theta_0\not\in\Lambda_{\tau}(\Gamma)$. Fix $\gamma_1$ minimal in
$\Gamma$ such that $\theta_0$ and $\gamma_1$ are non-comparable.
Let $\xi=\gamma_1\w\theta_0,\theta_1=\gamma_1\v\theta_0$. We
divide the proof into the following two cases.

\ni\textbf{Case 1:} $\gamma_1$ is the minimal element of $\Gamma$
(note that $\gamma_1$ is the minimal element of $D_\tau$ also).

We have $\xi\in \Lambda_{\tau}(\Gamma)$ (since $\xi<\gamma_1$);
further, $\xi\not\in D_{\tau}$ (again, since $\xi<\gamma_1$, the
minimal element of $D_\tau$).

\ni\textbf{Subcase 1(a):} Let $\theta_1$ be in $D_{\tau}$. Then
considering the diamond with $(\theta_1,\gamma_1),(\theta_0,\xi)$
as opposite sides, we have, $\theta_0\sim\xi$, and the result
follows (note that $\xi\in \Lambda_{\tau}(\Gamma)$).

\ni\textbf{Subcase 1(b):} Let $\theta_1\not\in D_{\tau}$. This
implies $\theta_0\in E_{\tau}$, and the result follows.

\ni\textbf{Case 2:} $\gamma_1$ is not the minimal element of
$\Gamma$.

 Let $\gamma_0\in\Gamma $ be such that $\gamma_0$ is covered by
$\gamma_1$ in $D_{\tau}$. Then $\xi\ge \gamma_0$ (since, both
$\theta_0$ and $\gamma_1$ are $>\gamma_0$).

\ni\textbf{Subcase 2(a):} Let $\xi>\gamma_0$. Then the fact that
$(\gamma_1,\gamma_0)$ is a cover in $D_\tau$ together with the
relation $\gamma_0<\xi<\gamma_1$ implies that $\xi\in
\Lambda_{\tau}(\Gamma)\,\setminus\,\Gamma$; in particular,
$\xi\not\in D_{\tau}$. Then as in Case 1, we obtain that
$\theta_0\in E_{\tau}$ if $\theta_1\not\in D_{\tau}$, and
$\theta_0\sim\xi$, if $\theta_1\in D_{\tau}$; and the result
follows (again note that $\xi\in \Lambda_{\tau}(\Gamma)$).

\ni\textbf{Subcase 2(b):} Let $\xi= \gamma_0$. We first note that
$\theta_1\not\in D_\tau$; for $\theta_1\in D_\tau$ would imply
that $\theta_0\in D_\tau$ (since $D_\tau$ is an embedded
sublattice (cf. Lemma \ref{embed})). Hence by induction
hypothesis, we obtain that
$$\theta_1\sim\eta, \mathrm{\ for\ some\ }\eta\in\Lambda_{\tau}(\Gamma)\cup E_\tau\leqno{(*)}$$
Also, considering the diamond with
$(\theta_1,\theta_0),(\gamma_1,\gamma_0)$ as opposite sides, we
have, $\theta_0\sim\theta_1$. Hence, the result follows in view of
(*).

 Note that the above proof includes the proof of the
starting point of induction, namely, $\theta_0$ is a maximal
element in ${\L}'_\tau$. To make this more precise, let
$\gamma_1,\xi,\theta_1$ be as above. Proceeding as above, we
obtain (by maximality of $\theta_0$) that $\theta_1\in D_\tau$.
Hence Subcases 1(b) and 2(b) do not exist (since in these cases
$\theta_1\not\in D_\tau$). In Subcases 1(a) and 2(a), we have
$\xi\in\Lambda_{\tau}(\Gamma)$, and $\theta_0\sim\xi$. The result
now follows in this case.
\end{pf}

Let
$Y_{\tau}(\Gamma)=\Lambda_{\tau}(\Gamma)\,\setminus\,E_{\tau}$.
Let us write
$$\begin{gathered}Y_{\tau}(\Gamma)=\Gamma\dot\cup Z_{\tau}(\Gamma)\\
G_{\tau}(\Gamma)=\{[\theta],\,\theta\in
Z_\tau(\Gamma),\,[\theta]\mathrm{\ is\ a\ non\
}E_\tau\mathrm{-class}\}
\end{gathered}$$
Combining the above Proposition with Lemmas \ref{two},\ref{one},
we obtain
\begin{prop}\label{five} $M_{\tau}/M^2_{\tau}$ is generated by
$\{{\overline{F}}_{[\theta]},[\theta]\in
G_{\tau}(\Gamma)\}\dot\cup\{{\overline{F}}_\gamma,\gamma\in\Gamma\}$.
\end{prop}
\section{A basis for the cotangent space at $P_{\tau}$}\label{cotgt} In this section,
we shall show that $\{{\overline{F}}_{[\theta]},[\theta]\in
G_{\tau}(\Gamma)\}\dot\cup\{{\overline{F}}_\gamma,\gamma\in\Gamma\}$
is in fact a basis for $M_{\tau}/M^2_{\tau}$.  We first recall
some basic facts on tangent cones and tangent spaces.

Let $X=Spec\,R\hookrightarrow {\mathbb{A}}^l$ be an affine
variety; let $S$ be the polynomial algebra $K[X_1,\cdots,X_l]$.
Let $P\in X$, and let $M_P$ be the maximal ideal in $K[X]$
corresponding to $P$ (we are concerned only with closed points of
$X$). Let $A={\mathcal{O}}_{X,P}$, the stalk at $P$; denote the
unique maximal ideal in $A$ by $M (=M_PR_{M_P})$. Then
$Spec\,gr(A,M)$, where
$gr(A,M)={\underset{j\in{\mathbb{Z}}_+}{\oplus}}\,M^j/M^{j+1}$ is
the {\emph{tangent cone}} to $X$ at $P$, and is denoted $TC_PX$.
Note that $$gr(R,M_P)=gr(A,M)$$

\ni\textbf{Tangent cone \& tangent space at $P$:} Let $I(X)$ be
the vanishing ideal of $X$ for the embedding
$X=Spec\,R\hookrightarrow {\mathbb{A}}^l$. Expanding $F\in I(X)$
in terms of the local co-ordinates at $P$, we have the following:

$\bullet$ $T_P(X)$, the tangent space to $X$ at $P$ is the zero
locus of the linear forms of $F$, for all $F\in I(X)$, i.e, the
degree one part in the (polynomial) local expression for $F$.

$\bullet$ $TC_P(X)$, the tangent cone to $X$ at $P$ is the zero
locus of the initial forms (i.e., form of smallest degree) of $F$,
for all $F\in I(X)$.

In the sequel, we shall denote the initial form of $F$ by $IN(F)$.

We collect below some well-known facts on singularities of $X$:

\vs.2cm\ni\textbf{Facts:} 1. $dim\,T_PX\ge dim\,X$ with equality
if and only if $X$ is smooth at $P$.


\ni 2. $X$ is smooth at $P$ if and only if $gr(R,M_P)$ is a
polynomial algebra.

\subsection{Determination of the degree one part of $J(\tau)$:} We
now take $X=X_{{\L}},P=P_\tau$; we shall denote $I=I(X_{{\L}}),
M_P=M_\tau$. As above, for $F\in M_\tau$, let ${\overline{F}}$
denote the class of $F$ in $M_{\tau}/M^2_{\tau}$. Let $J(\tau)$ be
the kernel of the surjective map
$$f_\tau:K[X_\theta,\theta\in{\L}]\rightarrow gr(R,M_\tau),\ X_\theta\mapsto {\overline{F_\theta}}$$
For $r\in\mathbb{N}$, let $f^{(r)}$ be the restriction of $f$ to
the degree $r$ part of the polynomial algebra
$K[X_\theta,\theta\in{\L}]$. We are interested in
$$f^{(1)}_\tau:{\underset{\theta\in{\L}}{\oplus}}\,KX_\theta\rightarrow
M_{\tau}/M^2_{\tau}$$ We shall first describe a complement to the
kernel of $f^{(1)}_\tau$, and then deduce a basis for
$M_{\tau}/M^2_{\tau}$ (which will turn out to be the set
$\{{\overline{F}}_{[\theta]},[\theta]\in
G_{\tau}(\Gamma)\}\dot\cup\{{\overline{F}}_\gamma,\gamma\in\Gamma\}$).

Let $V_\tau$ be the span of $\{{\overline{F_\beta}},\beta\in
D_\tau\}$ (in $M_{\tau}/M^2_{\tau}$); by Lemmas \ref{two}, we have
that $V_\tau$ is spanned by
$\{{\overline{F}}_\gamma,\gamma\in\Gamma\}$. Let $W_\tau$ be the
span of $\{{\overline{F_\theta}},\theta\not\in D_\tau\}$ (in
$M_{\tau}/M^2_{\tau}$).
\begin{lem}\label{vtau}
$\{{\overline{F}}_\gamma,\gamma\in\Gamma\}$ is a basis for
$V_\tau$.
\end{lem}
\begin{pf}
The surjective map $\pi_\tau:M_\tau/M^{2}_\tau\rightarrow
M'_\tau/M{'}^{2}_\tau$ induces a surjection $V_{\t}\rightarrow
M'_\tau/M{'}^{2}_\tau$, and the result follows from Corollary
\ref{dtau}.
\end{pf}
\begin{lem}\label{six}
The sum $V_\tau + W_\tau$ is direct.
\end{lem}
\begin{pf}
Let $v\in V_{\t}$; by Lemma \ref{vtau}, we can write
$v={\underset{\gamma\in\Gamma}{\sum}}\,a_\gamma
{\overline{F_\gamma}}$. Now if $v\in W_\tau$, then under the map
$$\pi_\tau:M_{\tau}/M^2_{\tau}\rightarrow M'_{\tau}/M{'}^2_{\tau}$$ we
obtain that $\pi_\tau(v)=0$ (since, $W_\tau\subseteq
ker\,\pi_\tau$). Hence we obtain that
$\pi_\tau({\underset{\gamma\in\Gamma}{\sum}}\,a_\gamma
{\overline{F_\gamma}})=0$ in $M'_{\tau}/M{'}^2_{\tau}$; this in
turn implies that $a_\gamma=0,\forall \gamma$ (cf. Corollary
\ref{dtau}). Hence $v=0$ and the required result follows.
\end{pf}
Let us write ${\underset{\theta\in{\L}}{\oplus}}\,KX_\theta
=A_\tau\oplus B_\tau$ where
$$A_\tau={\underset{\beta\in{D_\tau}}{\oplus}}\,KX_\beta,\ B_\tau={\underset{\theta\not\in{D_\tau}}{\oplus}}\,KX_\theta$$
Write $f^{(1)}_\tau=g^{(1)}_\tau +h^{(1)}_\tau$, where
$g^{(1)}_\tau$ (respectively $h^{(1)}_\tau$) is the restriction of
$f^{(1)}_\tau$ to $A_\tau$ (respectively $B_\tau$). Note that we
have surjections $$ g^{(1)}_\tau:A_\tau\rightarrow V_\tau,\
h^{(1)}_\tau:B_\tau\rightarrow W_\tau
$$ As a consequence of Lemma \ref{six}, we get the following
\begin{cor}\label{six'}
$ker\,f^{(1)}_\tau=ker\,g^{(1)}_\tau\oplus ker\,h^{(1)}_\tau$
\end{cor}
\begin{pf}
The inclusion``$\supseteq$" is clear. Let then $v\in
ker\,f^{(1)}_\tau$; write $v=a+b$, where $a\in A_\tau,b\in
B_\tau$. Then denoting $a'=f^{(1)}_\tau(a),b'=f^{(1)}_\tau(b)$, we
have, $0=a'+b'$. Also, we have that $a'\in V_\tau,b'\in W_\tau$.
Hence in view of Lemma \ref{six}, we obtain $a'=0=b'$. This
implies that $a\in ker\,g^{(1)}_\tau,b\in ker\,h^{(1)}_\tau$, and
the result follows.
\end{pf}
\begin{lem}\label{seven}
The span of $\{{X_\gamma},\gamma\in \Gamma\}$ is a complement to
the kernel of $g^{(1)}_\tau$.
\end{lem}
 \begin{pf} We have, $g^{(1)}_\tau(X_\beta)={\overline{F_\beta}}$. By Lemma \ref{vtau}, we have that
$\{{\overline{F_\gamma}},\gamma\in\Gamma\}$ is a basis for
$V_\tau$. The result now follows.
\end{pf}
Let $\{\xi_1,\cdots,\xi_r\}$ be a complete set of representatives
for $G_{\tau}(\Gamma)$.
\begin{lem}\label{eight}
The span of $\{{X_{\xi_{1}}},\cdots,{X_{\xi_{r}}}\}$ is a
complement to the kernel of $h^{(1)}_\tau$.
\end{lem}
\begin{pf}
A typical element $F$ of $I(X_{\L})$ such that IN($\,F$) is in
$B_\tau$ is of the form $$F=F_1+F_2$$ where $F_1$ is a linear sum
of diamond relations arising from diamonds having precisely one
vertex in $D_\tau$, and $F_2$ is a linear sum of diamond relations
arising from diamonds having precisely one side in $D_\tau$. Note
that in a typical term in $F_1$, the linear term is of the form
$a_\alpha X_\alpha,a_\alpha\in K$ where $\alpha\in E_\tau$;
similarly, in a typical term in $F_2$, the linear term is of the
form $b_{\theta\delta} (X_\theta- X_\delta),\,\theta,\delta\not\in
D_\tau,\theta\sim\delta,b_{\theta\delta}\in K$. Hence the kernel
of $h^{(1)}_\tau$ is generated by $\{X_\alpha,\alpha\in
E_\tau\}\cup\{(X_\theta- X_\delta),\theta\sim\delta\}$. The
required result now follows in view of Proposition \ref{four}.
\end{pf}
\begin{thm}\label{main2}
$\{{\overline{F}}_{[\theta]},[\theta]\in
G_{\tau}(\Gamma)\}\dot\cup\{{\overline{F}}_\gamma,\gamma\in\Gamma\}$
is a basis for $M_{\tau}/M^2_{\tau}$.
\end{thm}

\begin{pf} In view of Lemmas \ref{seven}, \ref{eight} and Corollary \ref{six'}, we obtain that

\ni
${\underset{\gamma\in\Gamma}{\oplus}}KX_\gamma\,\oplus\,{\underset{1\le
i\le r}{\oplus}}KX_{\xi_{i}}$ is a complement to the kernel of the
surjective map

\ni
$f^{(1)}_\tau:{\underset{\theta\in{\L}}{\oplus}}\,KX_\theta\rightarrow
M_{\tau}/M^2_{\tau}$. The result now follows.
\end{pf}



Let $T_{\tau}X_{\L}$ denote the tangent space to $X_{\L}$ at
$P_{\tau}$. As an immediate consequence of the above Theorem, we
have the following
\begin{cor}
$dim\,T_{\tau}X_{\L}=\#G_{\tau}(\Gamma)+\#\Gamma$.
\end{cor}

In view of the above Corollary, we obtain that $X_{\L}$ is
singular along $O_\tau$ if and only if
$\#G_{\tau}(\Gamma)+\#\Gamma>\#{\L}$. Let
$${\mathcal{S_{\L}}}=\{\tau<\sigma\,|\,\#G_{\tau}(\Gamma)+\#\Gamma>\#{\L}\}$$
(here, $\sigma$ is the cone associated to  $X_{\L}$). We obtain
from Proposition \ref{close} the following
\begin{thm}\label{main3}
Sing$\,X_{\L}={\underset{\tau\in\mathcal{S_{\L}}}{\cup}}\,X(D_\tau)$.
\end{thm}


\end{document}